\documentclass[12pt]{article}

\usepackage{a4wide}
\usepackage{amsfonts,amsmath}
\usepackage[utf8]{inputenc}

\newtheorem{theorem}{Theorem}[section]

\newtheorem{corollary}[theorem]{Corollary}
\newtheorem{lemma}[theorem]{Lemma}

\newtheorem{remark}[theorem]{Remark}
\newcommand{\proof} [1]
   { \noindent {\bf Proof.} #1 \hfill\rule{0.5em}{1.2ex} \par\medskip}

\newcommand{\R}{\mathbb{R}}

\newcommand{\N}{\mathbb{N}}
\newcommand{\norm}[1]{{\left\lVert{#1}\right\rVert}}
\newcommand{\spf}[2]{{\left\langle{#1},{#2}\right\rangle}}

\numberwithin{equation}{section}

\begin{document}

\title{A generalized inf--sup stable variational formulation \\ 
for the wave equation}
\author{Olaf~Steinbach$^1$, Marco~Zank$^2$}
\date{$^1$Institut f\"ur Angewandte Mathematik, TU Graz, \\ 
Steyrergasse 30, 8010 Graz, Austria \\[1mm]
{\tt o.steinbach@tugraz.at} \\[2mm]
$^2$Fakult\"at f\"ur Mathematik, Universit\"at Wien, \\
Oskar--Morgenstern--Platz 1, 1090 Wien, Austria \\[1mm]
{\tt marco.zank@univie.ac.at}}

\maketitle

\begin{abstract}
In this paper, we consider a variational formulation for the Dirichlet
problem of the wave equation with zero boundary and initial conditions,
where we use integration by parts in space and time. To prove unique
solvability in a subspace of $H^1(Q$) with $Q$ being the space--time
domain, the classical assumption is to consider the right--hand side
$f$ in $L^2(Q)$. Here, we analyze a generalized setting of this
variational formulation, which allows us to prove unique solvability
also for $f$ being in the dual space of the test space, i.e., the
solution operator is an isomorphism between the ansatz space and the
dual of the test space. This new approach is based on a suitable extension
of the ansatz space to include the information of the differential
operator of the wave equation at the initial time $t=0$. These results
are of utmost importance for the formulation and numerical analysis of
unconditionally stable space--time finite element methods, and for the
numerical analysis of boundary element methods to overcome the
well--known norm gap in the analysis of boundary integral operators.
\end{abstract}

\section{Introduction}
For the analysis of hyperbolic partial differential equations, a variety
of approaches such as Fourier methods, semigroups, or Galerkin methods
are available, see, e.g., \cite{Ladyzhenskaya1985, LionsMagenes1972,
Pazy1983, Wloka1982, Zeidler1990}. The theoretical results on existence
and uniqueness of solutions also form the basis for the numerical analysis
of related discretization schemes such as finite element methods, e.g.,
\cite{Bangerth2010, Bales1994, Cohen2002, DoerflerFindeisenWieners2016,
French1993, French1996, Koecher2014, Perugia2018, Richter1994,
SteinbachZankFEM2019, SteinbachZank2020, Zank2021, Zlotnik1994}
or boundary element methods, e.g.,
\cite{Aimi2009, Bamberger1986, Gimperlein2018,Joly2017,Sayas2016}.

As for elliptic second--order partial differential equations, we consider
the weak solution of the inhomogeneous wave equation in the energy space
$H^1(Q)$ with respect to the space--time domain $Q:=\Omega \times (0,T)$.
However, to ensure existence and uniqueness of a weak solution, we need
to assume $f \in L^2(Q)$, see, e.g., \cite[Theorem 5.1]{SteinbachZank2020}.
While this is a standard assumption to conclude sufficient regularity for
the solution $u$, and therefore, to obtain linear convergence for piecewise
linear finite element approximations $u_h$, stability of common finite
element discretizations require in most cases some CFL condition, which
relates the spatial and temporal mesh sizes to each other, see, e.g.,
\cite{SteinbachZankFEM2019, SteinbachZank2020, Zank2021, Zlotnik1994}.

Although the variational formulation to find $u$ in a suitable subspace
of $H^1(Q)$ is well--defined for $f$ being in the dual of the test space,
this is not sufficient to establish unique solvability. This is mainly
due to the missing information about the differential operator of the
wave equation at $t=0$ in the standard ansatz space. Hence, we are not
able to define an isomorphism between the ansatz space and the dual of
the test space. But such an isomorphism is an important ingredient in
the analysis of equivalent boundary integral formulations for boundary
value problems for the wave equation, and the numerical analysis
of related boundary and finite element methods.

In this paper, we are interested in inf--sup stable variational formulations
for the Dirichlet boundary value problem for the wave equation,
\begin{equation}\label{Welle}
    \left.
    \begin{array}{rcll}
        \Box u(x,t) & = & f(x,t) & 
        \text{for } (x,t) \in Q := \Omega \times (0,T), \\[1mm]
      u(x,t) & = & 0 & \text{for } (x,t)
                       \in \Sigma := \Gamma \times (0,T), \\[1mm]
        u(x,0) & = & 0 & \text{for } x \in \Omega, \\[1mm]
        \partial_t u(x,t)_{|t=0} & = & 0 & \text{for } x \in \Omega,
    \end{array}
    \right\}
\end{equation}
where $\Box u  := \partial_{tt}u - \Delta_x u$ is the wave operator,
$\Omega  \subset \R^d$, $d=1,2,3$, is a bounded domain with, for $d=2,3$,
Lipschitz boundary $\Gamma=\partial\Omega$, $T>0$ is a finite time horizon,
and $f$ is some given source. For simplicity, we assume homogeneous
Dirichlet boundary conditions as well as homogeneous initial conditions,
see, e.g., \cite{Ladyzhenskaya1985, Zank2021, Zlotnik1994} for the
treatment of inhomogeneous initial conditions. A possible variational
formulation of (\ref{Welle}) is to find  $u \in H^{1,1}_{0;0,}(Q) :=
L^2(0,T;H^1_0(\Omega)) \cap H^1_{0,}(0,T;L^2(\Omega))$ such that
\begin{equation} \label{Welle:VF}
    a(u,v) = \int_0^T \int_\Omega f(x,t) \, v(x,t) \, \mathrm dx \, \mathrm dt
\end{equation}
is satisfied for all $v \in H^{1,1}_{0;,0}(Q) := L^2(0,T;H^1_0(\Omega))
\cap H^1_{,0}(0,T;L^2(\Omega))$ with the bilinear form
$a(\cdot, \cdot) \colon \, H^{1,1}_{0;0,}(Q) \times H^{1,1}_{0;,0}(Q) \to \R$,
\begin{equation} \label{Welle:a}
  a(u,v) := - \int_0^T \int_\Omega \partial_t u(x,t) \, \partial_t v(x,t)
  \, \mathrm dx \, \mathrm dt +
  \int_0^T \int_\Omega \nabla_x u(x,t) \cdot \nabla_x v(x,t) \,
  \mathrm dx \, \mathrm dt
\end{equation}
for $u \in H^{1,1}_{0;0,}(Q), v \in H^{1,1}_{0;,0}(Q)$. In addition to the
standard Bochner space $L^2(0,T;H^1_0(\Omega))$, we use the space
$H^1_{0,}(0,T;L^2(\Omega))$ of all $v \in L^2(Q)$ with
$\partial_t v \in L^2(Q)$, and $v(x,0)=0$ for $x \in \Omega$.
Moreover, $H^1_{,0}(0,T;L^2(\Omega))$ is defined analogously
with $v(x,T)=0$ for $x \in \Omega$. The spaces $H^{1,1}_{0;0,}(Q)$,
$H^{1,1}_{0;,0}(Q)$ are Hilbert spaces with the inner products
\begin{equation*}
  \langle w, z \rangle_{H^{1,1}_{0;0,}(Q)} :=
  \langle w, z \rangle_{H^{1,1}_{0;,0}(Q)} :=
  \langle \partial_t w, \partial_t z \rangle_{L^2(Q)} +
  \langle \nabla_x w, \nabla_x z \rangle_{L^2(Q)},
\end{equation*}
and the corresponding induced norms $\| \cdot \|_{H^{1,1}_{0;0,}(Q)}$,
$\| \cdot \|_{H^{1,1}_{0;,0}(Q)}$. The bilinear form $a(\cdot , \cdot)$ in
\eqref{Welle:a} is continuous, i.e., for all
$u \in H^{1,1}_{0;0,}(Q)$ and $v \in H^{1,1}_{0;,0}(Q)$ we have
\begin{eqnarray*}
    |a(u,v)| &\leq&
    \| \partial_t u \|_{L^2(Q)} \| \partial_t v \|_{L^2(Q)}
    +
    \| \nabla_x u \|_{L^2(Q)} \| \nabla_x v \|_{L^2(Q)}  \nonumber \\[1mm]
    &\leq&  \sqrt{ \| \partial_t u \|^2_{L^2(Q)} + \| \nabla_x u \|_{L^2(Q)}^2}
           \sqrt{ \| \partial_t v \|^2_{L^2(Q)} + \| \nabla_x v \|_{L^2(Q)}^2}
           \nonumber \\[1mm]
             &=& \| u \|_{H^{1,1}_{0;0,}(Q)} \| v \|_{H^{1,1}_{0;,0}(Q)} .
                 \label{Welle:a:stetig}
\end{eqnarray*}
Note that the first initial condition $u(\cdot,0) = 0$ in $\Omega$ is
incorporated in the ansatz space $H^{1,1}_{0;0,}(Q)$, whereas the second
initial condition $\partial_t u(\cdot,t)_{|t=0} = 0$ in $\Omega$ is
considered as a natural condition in the variational formulation.
Thus, an inhomogeneous condition
$\mbox{$\partial_t u(\cdot,t)_{|t=0} = v_0$}$ in $\Omega$ with given $v_0$ could
be realized with the right--hand side $f_{v_0} \in [H^{1,1}_{0;,0}(Q)]'$,
\begin{equation} \label{Welle:fv0}
  \langle f_{v_0}, v \rangle_{Q} =
  \langle v_0, v(\cdot,0)\rangle_\Omega, \quad v \in H^{1,1}_{0;,0}(Q).
\end{equation}
However, known existence results for the variational formulation
\eqref{Welle:VF} do not allow right--hand sides in $[H^{1,1}_{0;,0}(Q)]'$,
which is the aim of the new approach as given in Section~\ref{Sec:Welle}. So,
when assuming $f \in L^2(Q)$, we are able to construct a unique solution
$u \in H^{1,1}_{0;0,}(Q)$ of the variational formulation \eqref{Welle:VF},
satisfying the stability estimate \cite[Theorem 5.1]{SteinbachZank2020},
see also \cite{Ladyzhenskaya1985, SteinbachZankFEM2019, Zank2020},
\begin{equation}\label{Welle:H1Stabilitaet}
  \| u \|_{H^{1,1}_{0;0,}(Q)} = \sqrt{\|\partial_tu\|^2_{L^2(Q)} +
  \| \nabla_x u \|^2_{L^2(Q)}}
  \leq \frac{1}{\sqrt{2}} \, T \, \| f \|_{L^2(Q)} \, .
\end{equation}
The stability estimate \eqref{Welle:H1Stabilitaet} does not fit to the
situation of the Banach--Ne\v{c}as--Babu\v{s}ka theorem as stated, e.g.,
in \cite[Theorem 2.6]{ErnGuermond2004}, see also \cite{Babuska1971,
BabuskaAziz1972, Necas1962}. The next theorem states that it is not
possible to prove the corresponding inf--sup condition, i.e.,
\eqref{Welle:InfSupH1}, for the bilinear form \eqref{Welle:a}.

\begin{theorem}{\rm \cite[Theorem 4.2.24]{Zank2020}}
  \label{Welle:Thm:InfSupH1}
  There does not exist a constant $c>0$ such that each right--hand side
  $f \in L^2(Q)$ and the corresponding solution $u \in H^{1,1}_{0;0,}(Q)$
  of \eqref{Welle:VF} satisfy
  \begin{equation*}
    \| u \|_{H^{1,1}_{0;0,}(Q)} \, \leq \, c \, \| f \|_{[H^{1,1}_{0;\,,0}(Q)]'} .
  \end{equation*}
  In particular, the inf--sup condition
  \begin{equation}\label{Welle:InfSupH1}
    c_S \, \| u \|_{H^{1,1}_{0;0,}(Q)} \, \leq \,
    \sup_{0 \neq v \in H^{1,1}_{0;\,,0}(Q)}
    \frac{| a(u,v) |}{\| v \|_{H^{1,1}_{0;,0}(Q)}} \quad
    \mbox{for all } u \in H^{1,1}_{0;0,}(Q)
    \end{equation}
    with a constant $c_S>0$ does not hold true.
\end{theorem}

\noindent
The proofs of the stability estimate \eqref{Welle:H1Stabilitaet} and
Theorem~\ref{Welle:Thm:InfSupH1} are based on an appropriate Fourier
analysis, using the eigenfunctions of the spatial differential operator
$-\Delta_x$, and the analysis of the related ordinary differential equation
\eqref{ODE}, which allows us also to present the essential ingredients for
the new approach. So, for $\mu >0$, we consider the scalar ordinary
differential equation
\begin{equation}\label{ODE}
    \Box_\mu u(t) :=
    \partial_{tt} u(t) + \mu u(t) = f(t) \quad \text{for } t \in (0,T), \quad
    u(0) = \partial_tu(t)_{|t=0} = 0 .
\end{equation}
The related variational formulation is to find 
$u \in H^1_{0,}(0,T)$ for a given right--hand side $f \in [H^1_{,0}(0,T)]'$
such that
\begin{equation}\label{ODE:VF}
 a_\mu(u,v) = \langle f, v \rangle_{(0,T)}
\end{equation}
is satisfied for all $v \in H^1_{,0}(0,T)$. The bilinear form
$a_\mu(\cdot, \cdot) \colon H^1_{0,}(0,T) \times H^1_{,0}(0,T) \to \R$
is defined by
\begin{equation}\label{ODE:a}
  a_\mu(u,v) := - \int_0^T \partial_t u(t) \, \partial_t v(t) \, \mathrm dt +
  \mu \int_0^T u(t) \, v(t) \, \mathrm dt
\end{equation}
for $u \in H^1_{0,}(0,T)$, $v \in H^1_{,0}(0,T)$. As before,
$u \in H^1_{0,}(0,T)$ covers the initial condition $u(0)=0$, while
$v \in H^1_{,0}(0,T)$ satisfies the terminal condition $v(T)=0$, and
where the inner product
$\langle \partial_t (\cdot), \partial_t (\cdot) \rangle_{L^2(0,T)}$
makes both to Hilbert spaces. Note that the second initial condition
$\partial_t u(t)_{|t=0}=0$ enters the variational formulation
\eqref{ODE:VF} as natural condition.
The dual space $[H^1_{,0}(0,T)]'$ is characterized as completion 
of $L^2(0,T)$ with respect to the Hilbertian norm
\begin{equation*}
    \| f \|_{[H^1_{,0}(0,T)]'} =
    \sup\limits_{0 \neq v \in H^1_{,0}(0,T)} \frac{|\langle f, v \rangle_{(0,T)}|}
    {\| \partial_t v \|_{L^2(0,T)}},
\end{equation*}
where $\langle \cdot , \cdot \rangle_{(0,T)}$ denotes the duality pairing
as extension of the inner product in $L^2(0,T),$ see, e.g.,
\cite[Satz 17.3]{Wloka1982}.
The continuity of the bilinear form \eqref{ODE:a} follows from
\begin{equation}
  |a_\mu(u,v)| \, \leq \,
  \left( 1 + \frac{4}{\pi^2} \mu T^2 \right) \,
  \| \partial_t u \|_{L^2(0,T)} \| \partial_t  v \|_{L^2(0,T)}
  \label{ODE:a:stetig}
\end{equation}
for all $u \in H^1_{0,}(0,T)$ and $v \in H^1_{,0}(0,T)$, where the
Cauchy--Schwarz and the Poincar\'e inequalities \cite[Lemma 3.4.5]{Zank2020}
are used. Furthermore, the bilinear form \eqref{ODE:a} satisfies the inf--sup
condition
\cite[Lemma 4.2]{SteinbachZank2020}
\begin{equation}\label{ODE:a_infsup}
    \frac{2}{2+\sqrt{\mu} T} \, \| \partial_t u \|_{L^2(0,T)} \leq
    \sup\limits_{0 \neq v \in H^1_{,0}(0,T)}
    \frac{| a_\mu(u,v) |}{\| \partial_t v \|_{L^2(0,T)}} \quad
    \text{for all } u \in H^1_{0,}(0,T) .
\end{equation}
Together with the positivity condition \cite[Lemma 4.2.4]{Zank2020}
\[
  a_\mu(u_v,v) = \langle v , v \rangle_{L^2(0,T)} =
  \| v \|^2_{L^2(0,T)} > 0
\]
for
\[
  u_v(t) = \frac{1}{\sqrt{\mu}} \int_0^t v(s) \sin (\sqrt{\mu}(t-s))
  \, \mathrm ds, \quad 0 \neq v \in L^2(0,T),
\]
we conclude unique solvability of \eqref{ODE:VF} as well as
the stability estimate
\begin{equation}\label{ODE stability}
\| \partial_t u \|_{L^2(0,T)} \leq \left( 1 + \frac{1}{2} \sqrt{\mu} T \right)
\| f \|_{[H^1_{,0}(0,T)]'}
\end{equation}
with the help of the Banach--Ne\v{c}as--Babu\v{s}ka theorem
\cite[Theorem 2.6]{ErnGuermond2004}. As discussed in
\cite[Remark 4.4]{SteinbachZank2020}, the stability estimate
(\ref{ODE stability}) is sharp in $\sqrt{\mu}$ and $T$, respectively. 
It turns out, however, that the estimate \eqref{ODE stability} is not
sufficient to prove a related stability estimate for the solution of
the wave equation \eqref{Welle}, see Theorem~\ref{Welle:Thm:InfSupH1}.
This is mainly due to the appearance of $\sqrt{\mu}$ in the stability
constant, i.e., (\ref{ODE stability}) is not uniform in $\mu$.
Instead, when assuming $f \in L^2(0,T)$, we can
prove the stability estimate \cite[Lemma 4.5]{SteinbachZank2020}
\begin{equation}\label{ODE stability L2}
\| \partial_t u \|^2_{L^2(0,T)} + \mu \, \| u \|^2_{L^2(0,T)} \leq 
\frac{1}{2} T^2 \, \| f \|^2_{L^2(0,T)} \, ,
\end{equation}
which allows to prove the stability estimate (\ref{Welle:H1Stabilitaet})
for the solution of the wave equation (\ref{Welle}), see
\cite[Theorem 5.1]{SteinbachZank2020} and \cite{Ladyzhenskaya1985, Zank2020}.

Since the variational formulation (\ref{ODE:VF}) is well--defined also
for $f \in [H^1_{,0}(0,T)]'$, we are interested to establish, instead of
(\ref{ODE:a_infsup}), an inf--sup stability condition with a constant,
which is independent of $\mu$, and which later on can be
generalized to the analysis of the variational problem
(\ref{Welle:VF}).

The remainder of this paper is structured as follows: In
Section~\ref{Sec:ODE}, we present the main ideas in order to solve the
ordinary differential equation \eqref{ODE}. For this purpose, we introduce a
suitable function space and prove a related inf--sup stability condition. Then,
in Section~\ref{Sec:Welle}, these results are generalized to analyze a
variational formulation for the wave equation \eqref{Welle}. The main result
of this paper is given in Theorem~\ref{Welle:Thm:ExistenzEindeutigkeit}, where
we state bijectivity of the solution operator for the Dirichlet problem of
the wave equation with zero boundary and initial conditions. Finally, in
Section~\ref{Sec:Zum}, we give some conclusions and comment on ongoing work.

\section{Variational formulation for the ODE} \label{Sec:ODE}
In this section, we derive a different setting of a variational
formulation for \eqref{ODE} in order to establish an inf--sup condition
with a constant independent of $\mu$. For this purpose, we first follow the
approach as for the heat equation, see, e.g., \cite{SteinbachZank2020}.
So, for given $u \in H^1_{0,}(0,T)$, we define
\begin{equation*}
  \langle \partial_{tt}u + \mu u, v \rangle_{(0,T)} :=
  a_\mu(u,v) \quad \text{for all } v \in H^1_{,0}(0,T).
\end{equation*}
Since the bilinear form $a_\mu (\cdot, \cdot)$ is continuous,
see \eqref{ODE:a:stetig}, $\partial_{tt}u + \mu u$ is a continuous functional
in $[H^1_{,0}(0,T)]'$. Therefore, by the Riesz representation theorem there
exists a unique element $w_u \in H^1_{,0}(0,T)$, satisfying
\begin{equation*}
  \langle \partial_{tt}u + \mu u, v \rangle_{(0,T)} =
  \langle \partial_t w_u, \partial_t v \rangle_{L^2(0,T)} \quad
  \text{for all } v \in H^1_{,0}(0,T),
\end{equation*}
and
\begin{equation} \label{ODE:NormH1Komma0}
  \| \partial_{tt}u + \mu u \|^2_{[H^1_{,0}(0,T)]'} =
  a_\mu(u,w_u) = \| \partial_t w_u \|^2_{L^2(0,T)}.
\end{equation}
At first glance, \eqref{ODE:NormH1Komma0} implies the inf--sup condition
\[
    \| \partial_{tt} u + \mu u \|_{[H^1_{,0}(0,T)]'}
    = \sup\limits_{0 \neq v \in H^1_{,0}(0,T)} 
    \frac{|a_\mu(u,v)|}{\|  \partial_t v \|_{L^2(0,T)}} \quad
    \text{ for } u \in H^1_{0,}(0,T),
\]
but $u \mapsto \| \partial_{tt} u + \mu u \|_{[H^1_{,0}(0,T)]'}$ only defines a 
semi--norm in $H^1_{0,}(0,T)$ since the differential operator
$\partial_{tt} + \mu$ is treated only in $(0,T)$, i.e., its behavior in
$t=0$ is not covered in \eqref{ODE:NormH1Komma0}. As example, consider,
e.g., $u \in H^1_{0,}(0,T)$ with $u(t) = \sin (\sqrt{\mu} t)$ and
$\| \partial_{tt} u + \mu u \|_{[H^1_{,0}(0,T)]'} = 0$. Hence, we need to
modify the ansatz space to determine $u$ in a suitable way. For this
purpose, we first introduce notations for additional function spaces and
operators.

In this work, $C_0^\infty(A)$ is the set of infinitely differentiable
real--valued functions with compact support in any domain $A \subset \R^d$,
$d=1,2,3,4.$ The set $C_0^\infty(A)$ is endowed with the, usual for
distributions, locally convex topology and is called the space of test
functions on $A.$ The set of (Schwartz) distributions $[C_0^\infty(A)]'$
is given by all linear and sequentially continuous functionals on
$C_0^\infty(A)$.

For given $u \in L^2(0,T)$, we define the extension
$\widetilde{u} \in L^2(-T,T)$ by
\begin{equation}\label{ODE:DefNullFortsetzung}
  \widetilde{u}(t) \, := \, \begin{cases}
                            u(t) & \text{for } t \in (0,T), \\
                               0 & \text{for } t \in (-T,0].
                      \end{cases}
\end{equation}
The application of the differential operator $\Box_\mu$
to $\widetilde{u}$ is defined as distribution on $(-T,T)$, i.e., for 
all test functions $\varphi \in C^\infty_0(-T,T)$, we define
\begin{equation}\label{ODE:Def_verallgemeinerte_Ableitung}
  \langle \Box_\mu \widetilde{u} , \varphi \rangle_{(-T,T)} :=
  \int_{-T}^T \widetilde{u}(t) \, \Box_\mu \varphi(t) \, \mathrm dt =
  \int_0^T u(t) \, \Box_\mu \varphi(t) \, \mathrm dt.
\end{equation}
This motivates to consider the dual space $[H^1_0(-T,T)]'$ of
$H^1_0(-T,T)$, which is characterized as completion of $L^2(-T,T)$
with respect to the Hilbertian norm
\begin{equation*}
  \| g \|_{[H^1_0(-T,T)]'} :=
  \sup_{0 \neq z \in H^1_0(-T,T) }
  \frac{|\langle g, z \rangle_{(-T,T)}|}{ \| \partial_t z \|_{L^2(-T,T)}},
\end{equation*}
where $\langle \cdot , \cdot \rangle_{(-T,T)}$ denotes the duality pairing
as extension of the inner product in $L^2(-T,T)$, see, e.g.,
\cite[Satz 17.3]{Wloka1982}. In other words, for $[H^1_0(-T,T)]'$, there
exists an inner product $\spf{\cdot}{\cdot}_{[H^1_0(-T,T)]'},$ inducing the
norm $\| \cdot \|_{[H^1_0(-T,T)]'} = \sqrt{\spf{\cdot}{\cdot}_{[H^1_0(-T,T)]'}}$,
i.e., with this abstract inner product, $[H^1_0(-T,T)]'$ is a Hilbert space.
Additionally, we define the subspace
\begin{multline*}
 H^{-1}_{[0,T]}(-T,T) := \Big\{ g \in [H^1_0(-T,T)]' \colon  \\ 
 \forall z \in H^1_0(-T,T) \text{ with } \mathrm{supp}\, z \subset (-T,0)
 \colon  \, \langle g, z \rangle_{(-T,T)} = 0 \Big\} \subset
 [H^1_0(-T,T)]',
\end{multline*}
endowed with the Hilbertian norm $\norm{\cdot}_{[H^1_0(-T,T)]'}.$ To
characterize the subspace $H^{-1}_{[0,T]}(-T,T)$, we introduce the
following notations. Let
${\mathcal{R}} : H^1_0(-T,T) \to H^1_{,0}(0,T)$ be the continuous
and surjective restriction operator, defined by $\mathcal R z = z_{|(0,T)}$
for $z \in H^1_0(-T,T)$, with its adjoint operator
$\mathcal R' : [H^1_{,0}(0,T)]' \to [H^1_0(-T,T)]'$. Furthermore,
let $\mathcal E : H^1_{,0}(0,T) \to H^1_0(-T,T)$ be any continuous
and injective extension operator with its adjoint operator
$\mathcal E' : [H^1_0(-T,T)]' \to [H^1_{,0}(0,T)]'$, satisfying
\begin{equation*}
  \| \mathcal E v \|_{H^1_0(-T,T)} \leq c_{\mathcal E} \| v \|_{H^1_{,0}(0,T)}
\end{equation*}
with a constant $c_{\mathcal E} > 0$ and $\mathcal R \mathcal E v = v$
for all $v \in H^1_{,0}(0,T).$ An example for such an extension operator is
given by reflection in $t=0$, i.e., consider the function $\overline{v}$,
defined for $v \in H^1_{,0}(0,T)$ by
\begin{equation*}
  \overline{v}(t) = \begin{cases}
                        v(t)  & \text{for } t \in [0,T), \\
                        v(-t) & \text{for } t \in (-T,0),
                   \end{cases}
\end{equation*}
which leads to the constant $c_{\mathcal E} = 2$ in this particular
case. With this, we prove the following lemma.

\begin{lemma} \label{ODE:Lem:DualH10Support}
  The spaces $(H^{-1}_{[0,T]}(-T,T), \| \cdot \|_{[H^1_0(-T,T)]'})$ and
  $([H^1_{,0}(0,T)]', \| \cdot \|_{[H^1_{,0}(0,T)]'})$ are isometric, i.e.,
  the mapping
  \begin{equation*}
    \mathcal E'_{|H^{-1}_{[0,T]}(-T,T)} \colon \,
    H^{-1}_{[0,T]}(-T,T) \to  [H^1_{,0}(0,T)]'
  \end{equation*}
  is bijective with
  \begin{equation*}
    \| g \|_{[H^1_0(-T,T)]'} = \| \mathcal E' g \|_{[H^1_{,0}(0,T)]'}
    \quad \text{ for all } g \in H^{-1}_{[0,T]}(-T,T).
  \end{equation*}
     In addition, for $g \in H^{-1}_{[0,T]}(-T,T)$, the relation
    \begin{equation} \label{ODE:Lem:DualH10Support:gDurchR}
      \langle g, z \rangle_{(-T,T)} =
      \langle \mathcal E' g, \mathcal R z \rangle_{(0,T)} \quad
      \mbox{for all } z \in H^1_0(-T,T)
    \end{equation}
    holds true, i.e.,
    $\mathcal R' \mathcal E' g = g$. In particular, the subspace
    $H^{-1}_{[0,T]}(-T,T) \subset [H^1_0(-T,T)]'$ is closed, i.e., complete.
\end{lemma}
  
\proof{ 
  First, we prove that
  $\| g \|_{[H^1_0(-T,T)]'} = \| \mathcal E' g \|_{[H^1_{,0}(0,T)]'}$
  for all $g \in H^{-1}_{[0,T]}(-T,T)$. For this purpose, let
  $g \in H^{-1}_{[0,T]}(-T,T)$ be arbitrary but fixed. The Riesz
  representation theorem gives the unique element $z_g \in H^1_0(-T,T)$ with
  \begin{equation*}
    \langle g, z \rangle_{(-T,T)} =
    \langle \partial_t z_g, \partial_t z \rangle_{L^2(-T,T)} \quad
    \text{ for all } z \in H^1_0(-T,T),
  \end{equation*}
  and $\| g \|_{[H^1_0(-T,T)]'} = \| \partial_t z_g \|_{L^2(-T,T)}$.
  It holds true that $z_{g|(-T,0)} = 0$, since we have
  \begin{equation*}
    0 = \langle g, z \rangle_{(-T,T)} =
    \langle \partial_t z_g, \partial_t z \rangle_{L^2(-T,T)} =
    \langle \partial_t z_{g|(-T,0)}, \partial_t z_{|(-T,0)} \rangle_{L^2(-T,0)}
  \end{equation*}
  for all $z \in H^1_0(-T,T)$ with $\mathrm{supp}\, z \subset (-T,0)$.
  Hence, we have
  \begin{equation}\label{ODE:Lem:DualH10Support:Beweis:gDurchR}
    \langle g, z \rangle_{(-T,T)} =
    \langle \partial_t z_g, \partial_t z \rangle_{L^2(-T,T)} =
    \langle \partial_t \mathcal R z_g,
    \partial_t \mathcal R z \rangle_{L^2(0,T)}
  \end{equation}
  for all $z \in H^1_0(-T,T)$. So, using
  \eqref{ODE:Lem:DualH10Support:Beweis:gDurchR} with $z=\mathcal E v$
  for $v \in H^1_{,0}(0,T)$ this gives
  \begin{equation} \label{ODE:Lem:DualH10Support:Beweis:E'gDurchR}
    \langle \mathcal E'g, v \rangle_{(0,T)} =
    \langle g, \mathcal E v \rangle_{(-T,T)} =
    \langle \partial_t \mathcal R z_g,
    \partial_t \mathcal R \mathcal E v \rangle_{L^2(0,T)} =
    \langle \partial_t \mathcal R z_g, \partial_t v \rangle_{L^2(0,T)} ,
  \end{equation}
  i.e.,
  \[
    \| {\mathcal{E}}' g \|_{[H^1_{,0}(0,T)]'} =
    \| \partial_t {\mathcal{R}} z_g \|_{L^2(0,T)} =
    \| \partial_t z_g \|_{L^2(-T,T)} =
    \| g \|_{[H^1_0(-T,T)]'} .
  \]
  Second, we prove that $\mathcal E'_{|H^{-1}_{[0,T]}(-T,T)} $ is
  surjective. For this purpose, let $f \in [H^1_{,0}(0,T)]'$ be given.
  Set $g_f = \mathcal R'f \in [H^1_0(-T,T)]'$, i.e.,
  \begin{equation*}
    \langle g_f, z \rangle_{(-T,T)} =
    \langle {\mathcal{R}}' f , z \rangle_{(-T,T)} =
    \langle f, \mathcal R z \rangle_{(0,T)} 
  \end{equation*}
  for all $z \in H^1_0(-T,T)$. With this it follows immediately that
  $g_f \in H^{-1}_{[0,T]}(-T,T)$. Moreover, we have
  \begin{equation*}
    \langle \mathcal E' g_f, v \rangle_{(0,T)} =
    \langle  g_f, \mathcal E v \rangle_{(-T,T)} =
    \langle f, \mathcal R \mathcal E v \rangle_{(0,T)} =
    \langle f, v \rangle_{(0,T)}
  \end{equation*}
  for all $v \in H^1_{,0}(0,T)$, i.e., $\mathcal E' g_f = f$ in
  $[H^1_{,0}(0,T)]'$. In other words,
  $\mathcal E'_{|H^{-1}_{[0,T]}(-T,T)} $ is surjective.
  Third, the equality \eqref{ODE:Lem:DualH10Support:gDurchR} follows
  from \eqref{ODE:Lem:DualH10Support:Beweis:gDurchR} and
  \eqref{ODE:Lem:DualH10Support:Beweis:E'gDurchR} for $v = \mathcal R z$
  for any $z \in H^1_0(-T,T)$.
  The last assertion of the lemma is straightforward.}

\noindent
The last lemma gives immediately the following corollary.

\begin{corollary} \label{ODE:Kor:NormHSupp0T}
  For all $g \in H^{-1}_{[0,T]}(-T,T)$, the norm representation
  \begin{equation*} 
    \| g \|_{[H^1_0(-T,T)]'} =
    \sup\limits_{0 \neq v \in H^1_{,0}(0,T)}
    \frac{| \langle g, \mathcal E v \rangle_{(-T,T)} |}
    {\| \partial_t v \|_{L^2(0,T)}}
  \end{equation*}
  holds true.
\end{corollary}

\proof{Let $g \in H^{-1}_{[0,T]}(-T,T)$ be arbitrary but fixed. With
  Lemma~\ref{ODE:Lem:DualH10Support}, we have
  \begin{equation*}
    \| g \|_{[H^1_0(-T,T)]'} =
    \| \mathcal E' g \|_{[H^1_{,0}(0,T)]'} =
    \sup\limits_{0 \neq v \in H^1_{,0}(0,T)}
    \frac{| \langle \mathcal E' g , v\rangle_{(0,T)} |}
    {\| \partial_t v \|_{L^2(0,T)}}  =
    \sup\limits_{0 \neq v \in H^1_{,0}(0,T)}
    \frac{| \langle g , \mathcal E v\rangle_{(-T,T)} |}
    {\| \partial_t v \|_{L^2(0,T)}},
  \end{equation*}
  i.e., the assertion is proven.}

\noindent
Next, we introduce
\[
  \mathcal H (0,T) := \Big \{ u= \widetilde{u}_{|(0,T)} : 
  \widetilde{u} \in L^2(-T,T), \; \widetilde{u}_{|(-T,0)} = 0, \;
  \Box_\mu \widetilde{u} \in [H^1_0(-T,T)]' \Big \}
\]
with the norm
\[
  \| u \|_{\mathcal{H}(0,T)} := \sqrt{\| u \|^2_{L^2(0,T)} + 
    \| \Box_\mu \widetilde{u} \|^2_{[H^1_0(-T,T)]'} } \; .
\]
For a function $u \in \mathcal{H}(0,T)$, the condition
$\Box_\mu \widetilde{u} \in [H^1_0(-T,T)]'$ involves that there
exists an element $f_u \in [H^1_0(-T,T)]'$ with
\begin{equation*}
  \langle \Box_\mu \widetilde{u} , \varphi \rangle_{(-T,T)} =
  \spf{f_u}{\varphi}_{(-T,T)} \quad \text{ for all }
  \varphi \in  C^\infty_0(-T,T).
\end{equation*}
Note that $\varphi \in H^1_0(-T,T)$ for $\varphi \in C^\infty_0(-T,T)$,
and that $C^\infty_0(-T,T)$ is dense in $H^1_0(-T,T)$. Hence, the
element $f_u \in [H^1_0(-T,T)]'$ is unique and therefore, in the
following, we identify the distribution
$\Box_\mu \widetilde{u} \colon \, C^\infty_0(-T,T) \to \R$ with the
functional $f_u \colon \, H^1_0(-T,T) \to \R$.

Next, we state properties of the space $\mathcal H (0,T)$. Clearly,
$(\mathcal H (0,T),\| \cdot \|_{\mathcal H (0,T)})$ is a normed
vector space, and it is even a Banach space.

\begin{lemma}\label{ODE:Lem:Banachraum}
The normed vector space $(\mathcal H (0,T),\| \cdot \|_{\mathcal H (0,T)})$
is a Banach space.
\end{lemma}

\proof{Consider a Cauchy sequence
  $(u_n)_{n \in {\N }} \subset \mathcal H (0,T)$. Hence,
  $(u_n)_{n \in {\N }} \subset L^2(0,T)$ is also a Cauchy sequence
  in $L^2(0,T)$, and $(\Box_\mu \widetilde{u}_n)_{n \in {\N }} \subset
  [H^1_0(-T,T)]'$ is also a Cauchy sequence in $[H^1_0(-T,T)]'$. So,
  there exist $u \in L^2(0,T)$ and $ f \in [H^1_0(-T,T)]'$ with
  \[
    \lim\limits_{n \to \infty} \| u_n - u \|_{L^2(0,T)} = 0, \quad
    \lim\limits_{n \to \infty}
    \| \Box_\mu \widetilde{u}_n - f \|_{[H^1_0(-T,T)]'} = 0.
  \]
  For $\varphi \in C_0^\infty(-T,T)$, we have
\begin{align*}
  \langle \Box_\mu \widetilde{u}, \varphi \rangle_{(-T,T)}
  &= \langle \widetilde{u}, \Box_\mu \varphi \rangle_{L^2(-T,T)} =
    \int_0^T u(t) \, \Box_\mu \varphi(t) \, \mathrm dt =
    \lim\limits_{n \to \infty}
    \int_0^T u_n(t) \, \Box_\mu \varphi(t) \, \mathrm dt \\ 
  &= \lim\limits_{n \to \infty}
    \langle \widetilde{u}_n , \Box_\mu \varphi \rangle_{L^2(-T,T)} =
    \lim\limits_{n \to \infty}
    \langle \Box_\mu \widetilde{u}_n , \varphi \rangle_{(-T,T)} \, = \,
    \langle f , \varphi \rangle_{(-T,T)},
\end{align*}
i.e., $\Box_\mu \widetilde{u} = f \in [H^1_0(-T,T)]'$. Hence,
$u \in \mathcal H (0,T)$ follows.}

\noindent
With the abstract inner product
$\langle \cdot, \cdot \rangle_{[H^1_0(-T,T)]'}$ of $[H^1_0(-T,T)]'$, the
inner product
\begin{equation*}
  \langle u, v \rangle_{\mathcal H (0,T)} :=
  \langle u, v \rangle_{L^2 (0,T)} +
  \langle \Box_\mu \widetilde{u}, \Box_\mu \widetilde{v}
  \rangle_{[H^1_0(-T,T)]'}, \quad u,v \in {\mathcal{H}}(0,T),
\end{equation*}
induces the norm $\| \cdot \|_{\mathcal H (0,T)}$. Hence, the space
$(\mathcal H (0,T), \langle \cdot, \cdot \rangle_{\mathcal H (0,T)})$ is
even a Hilbert space, but this abstract inner product is not used
explicitly in the remainder of this work.

\begin{lemma} \label{ODE:Lem:BoxInH10Support}
    For all $u \in \mathcal H (0,T)$ there holds
    $\Box_\mu \widetilde{u} \in H^{-1}_{[0,T]}(-T,T)$ and
    \begin{equation} \label{ODE:BoxNormdarstellung}
        \| \Box_\mu \widetilde{u} \|_{[H^1_0(-T,T)]'} =  \sup\limits_{0 \neq v
          \in H^1_{,0}(0,T)} \frac{| \langle \Box_\mu \widetilde{u}, \mathcal E
          v \rangle_{(-T,T)} |}{\| \partial_t v \|_{L^2(0,T)}} .
    \end{equation}
\end{lemma}
  
\proof{First, we prove that $\Box_\mu \widetilde{u} \in H^{-1}_{[0,T]}(-T,T)$.
  For this purpose, let $u \in \mathcal H (0,T)$ and $z \in H^1_0(-T,T)$
  with $\mathrm{supp}\, z \subset (-T,0)$ be arbitrary but fixed. Due to
  $z_{|(-T,0)} \in H^1_0(-T,0)$ there exists a sequence
  $(\psi_n)_{n \in \N} \subset C^\infty_0(-T,0)$ with
  $\|\partial_t z_{|(-T,0)} - \partial_t \psi_n \|_{L^2(-T,0)} \to 0$ as
  $n \to \infty$. For $n \in \N$, define
  \begin{equation*}
    \varphi_n(t) =
    \begin{cases}
        \psi_n(t) & \text{for } t \in (-T,0), \\
        0       & \text{for } t \in [0,T),
    \end{cases}
  \end{equation*}
  i.e., $(\varphi_n)_{n \in \N} \subset C^\infty_0(-T,T)$ satisfies
  \begin{equation*}
    \|\partial_t z - \partial_t \varphi_n \|_{L^2(-T,T)} =
    \|\partial_t z_{|(-T,0)} - \partial_t \psi_n \|_{L^2(-T,0)} \to 0
  \end{equation*}
  as $n \to \infty.$ So, it follows that
  \begin{equation*}
    \langle \Box_\mu \widetilde{u}, z \rangle_{(-T,T)} =
    \lim_{n\to\infty} \langle \Box_\mu \widetilde{u}, \varphi_n \rangle_{(-T,T)}
    = \lim_{n\to\infty} \int_0^T u(t) \, \Box_\mu \varphi_n(t) \, \mathrm dt = 0,
\end{equation*}
and therefore, the assertion. The norm representation follows
from $\Box_\mu \widetilde{u} \in H^{-1}_{[0,T]}(-T,T)$ and
Corollary~\ref{ODE:Kor:NormHSupp0T}.}

\begin{lemma}\label{ODE:Lem:Darstellung_Bilinearform}
  It holds true that $H^1_{0,}(0,T) \subset \mathcal H (0,T)$. Furthermore,
  each $u \in H^1_{0,}(0,T)$ with zero extension $\widetilde{u}$,
  as defined in \eqref{ODE:DefNullFortsetzung}, satisfies
  \begin{equation} \label{ODE:NormBoxuDurchAbleitungu}
    \| \Box_\mu \widetilde{u} \|_{[H^1_0(-T,T)]'} \leq
    \left( 1 + \frac{4}{\pi^2} \mu T^2 \right) \,
    \| \partial_t u \|_{L^2(0,T)} ,
  \end{equation}
  and
  \begin{equation}\label{ODE:Darstellung_Bilinearform}
    \langle \Box_\mu \widetilde{u} , z \rangle_{(-T,T)} =
    a_\mu(u,\mathcal R z) =
    - \langle \partial_t u , \partial_t \mathcal R z \rangle_{L^2(0,T)} +
    \mu \, \langle u , \mathcal R z \rangle_{L^2(0,T)}
  \end{equation}
  for all $z \in H^1_0(-T,T),$ where $a_\mu(\cdot,\cdot)$ is the
  bilinear form \eqref{ODE:a}.
\end{lemma}

\proof{First, we prove that $H^1_{0,}(0,T) \subset \mathcal H (0,T)$.
  For $u \in H^1_{0,}(0,T)$, we define the extension $\widetilde{u}$,
  see \eqref{ODE:DefNullFortsetzung}. By construction, we have
  $\widetilde{u} \in L^2(-T,T)$, and $\widetilde{u}_{|(-T,0)}=0$. It
  remains to prove that $\Box_\mu \widetilde{u} \in [H^1_0(-T,T)]'$.
  For this purpose, define the functional $f_u \in [H^1_0(-T,T)]'$ by
  \begin{equation*}
    \langle f_u, z \rangle_{(-T,T)} :=  a_\mu(u,\mathcal R z)
  \end{equation*}
  for all $z \in H^1_0(-T,T)$, where $a_\mu(\cdot,\cdot)$ is the bilinear
  form \eqref{ODE:a}. The continuity of $f_u$ follows from
  \begin{align*}
    |\langle f_u, z \rangle_{(-T,T)}| = | a_\mu(u,\mathcal R z) |
    & \leq \left( 1 + \frac{4}{\pi^2} \mu T^2 \right) \,
      \| \partial_t u \|_{L^2(0,T)}  \| \partial_t  \mathcal R z \|_{L^2(0,T)}
    \\
    & \leq \left( 1 + \frac{4}{\pi^2} \mu T^2 \right) \,
      \| \partial_t u \|_{L^2(0,T)} \| \partial_t z \|_{L^2(-T,T)}  
\end{align*}
for all $z \in H^1_0(-T,T)$, where the estimate \eqref{ODE:a:stetig} is
used. Using the definition (\ref{ODE:Def_verallgemeinerte_Ableitung}), and
integration by parts, this gives
\begin{align*}
  \langle \Box_\mu \widetilde{u} , \varphi \rangle_{(-T,T)}
  &= \int_0^T u(t)\, \Box_\mu \varphi(t) \, \mathrm dt =
  - \langle \partial_t u , \partial_t \mathcal R \varphi \rangle_{L^2(0,T)}
  + \mu \langle u ,  \mathcal R \varphi \rangle_{L^2(0,T)} \\
   &= \langle f_u, \varphi \rangle_{(-T,T)}
\end{align*}
for all $\varphi \in C^\infty_0(-T,T)$, i.e.,
$\Box_\mu \widetilde{u} = f_u \in [H^1_0(-T,T)]'.$ The equality
\eqref{ODE:Darstellung_Bilinearform} follows from the density of
$C^\infty_0(-T,T)$ in $H^1_0(-T,T)$. The estimate
\eqref{ODE:NormBoxuDurchAbleitungu} is proven by
\begin{align*}
  \| \Box_\mu \widetilde{u} \|_{[H^1_0(-T,T)]'}
  &= \sup\limits_{0 \neq v \in H^1_{,0}(0,T)}
    \frac{| \langle \Box_\mu \widetilde{u}, \mathcal E v \rangle_{(-T,T)}|}
    {\| \partial_t v \|_{L^2(0,T)}} =
    \sup\limits_{0 \neq v \in H^1_{,0}(0,T)}
    \frac{| \langle f_u, \mathcal E v \rangle_{(-T,T)} |}
    {\| \partial_t v \|_{L^2(0,T)}} \\
  &= \sup\limits_{0 \neq v \in H^1_{,0}(0,T)}
    \frac{ |a_\mu(u,\mathcal R \mathcal E v) | }
    {\| \partial_t v \|_{L^2(0,T)}} \, \leq \,
    \left( 1 + \frac{4}{\pi^2} \mu T^2 \right) \,
    \| \partial_t u \|_{L^2(0,T)} ,
\end{align*}
when using the norm representation \eqref{ODE:BoxNormdarstellung},
the equality \eqref{ODE:Darstellung_Bilinearform}, and the bound
\eqref{ODE:a:stetig}.}

\noindent
Next, by completion, we define the Hilbert space
\[
  \mathcal H_{0,}(0,T) :=
  \overline{H^1_{0,}(0,T)}^{\| \cdot \|_{\mathcal H (0,T)}}
  \subset \mathcal H (0,T),
\]
endowed with the Hilbertian norm $\norm{\cdot}_{\mathcal H(0,T)}$, i.e.,
\begin{equation*}
  \mathcal H_{0,}(0,T) = \Big \{ v \in \mathcal H(0,T) \colon \,
  \exists (v_n)_{n \in \N} \subset H^1_{0,}(0,T) \text{ with }
  \lim\limits_{n \to \infty }\norm{v_n - v}_{\mathcal H(0,T)} = 0 \Big \}.
\end{equation*}

\begin{lemma}\label{ODE:Lem:Friedrichs}
For $u \in \mathcal H_{0,}(0,T)$ there holds
\[
  \| \Box_\mu \widetilde{u} \|_{[H^1_0(-T,T)]'} \geq
  \frac{\sqrt{2}}{T} \, \| u \|_{L^2(0,T)} .
\]
\end{lemma}

\proof{For $0 \neq u \in \mathcal H_{0,}(0,T)$, there exists a
non--trivial sequence $(u_n)_{n \in \N } \subset H^1_{0,}(0,T)$,
$u_n \not\equiv 0$, with
\[
  \lim\limits_{n \to \infty} \| u - u_n \|_{\mathcal H (0,T)} = 0.
\]
For each $u_n \in H^1_{0,}(0,T)$, we define $w_n \in H^1_{,0}(0,T)$ as
unique solution of the backward problem
\begin{equation}\label{ODE backward}
  \partial_{tt} w_n(t) + \mu w_n(t) = u_n(t) \quad \text{for } t \in (0,T),
  \quad w_n(T)=\partial_tw_n(t)_{|t=T} =0,
\end{equation}
i.e., $w_n \in H^1_{,0}(0,T)$ solves the variational problem
\[
  a_\mu(v,w_n) = \langle u_n , v \rangle_{L^2(0,T)} \quad
  \text{for all } v \in H^1_{0,}(0,T)
\]
with the bilinear form \eqref{ODE:a}. In particular for $v=u_n$, this gives
\[
  a_\mu(u_n,w_n) = \| u_n \|^2_{L^2(0,T)} .
\]
Analogously to the estimate \eqref{ODE stability L2} for the solution of
\eqref{ODE}, we conclude
\[
  \| \partial_t w_n \|^2_{L^2(0,T)} + \mu \, \| w_n \|^2_{L^2(0,T)} \leq
  \frac{1}{2} \, T^2 \, \| u_n \|^2_{L^2(0,T)}
\]
for the solution $w_n$ of (\ref{ODE backward}), i.e.,
\[
  \| \partial_t w_n \|_{L^2(0,T)} \leq
  \frac{1}{\sqrt{2}} \, T \, \| u_n \|_{L^2(0,T)} .
\]
For the zero extension $\widetilde{u}_n \in L^2(-T,T)$ of
$u_n \in H^1_{0,}(0,T)$, we obtain, when using the norm representation
\eqref{ODE:BoxNormdarstellung}, and \eqref{ODE:Darstellung_Bilinearform},
that
\begin{align*}
  \| \Box_\mu \widetilde{u}_n \|_{[H^1_0(-T,T)]'}
  &= \sup\limits_{0 \neq v \in H^1_{,0}(0,T)}
    \frac{| \langle \Box_\mu \widetilde{u}_n, \mathcal E v \rangle_{(-T,T)}|}
    {\| \partial_t v \|_{L^2(0,T)}} \geq
    \frac{| \langle \Box_\mu \widetilde{u}_n , \mathcal E w_n \rangle_{(-T,T)}|}  
    {\| \partial_t w_n \|_{L^2(0,T)}} \\
  &= \frac{| a_\mu(u_n,w_n) |}{\| \partial_t w_n \|_{L^2(0,T)}} \, = \,
    \frac{\| u_n \|^2_{L^2(0,T)}}{\| \partial_t w_n \|_{L^2(0,T)}} \, \geq \,
    \frac{\sqrt{2}}{T} \, \| u_n \|_{L^2(0,T)} ,
\end{align*}
and the assertion follows by completion for $n \to \infty$.}

\begin{corollary} 
  The inner product space
  $\left(\mathcal H_{0,}(0,T), \langle \Box_\mu \widetilde{(\cdot)},
    \Box_\mu \widetilde{(\cdot)} \rangle_{[H^1_0(-T,T)]'} \right)$ is
  complete, i.e., a Hilbert space.
\end{corollary}

\proof{The assertion follows immediately from Lemma~\ref{ODE:Lem:Friedrichs}.}

\noindent
In the following, $\mathcal H_{0,}(0,T)$ is endowed with the Hilbertian
norm $\| \Box_\mu \widetilde{(\cdot)} \|_{[H^1_0(-T,T)]'}$. With this new
Hilbert space, the bilinear form
\[
  \widetilde{a}_\mu(\cdot,\cdot) \colon \,
  \mathcal H_{0,}(0,T) \times H^1_{,0}(0,T) \to \R, \quad 
  \widetilde{a}_\mu(u,v) := \langle \Box_\mu \widetilde{u} ,
  \mathcal E v \rangle_{(-T,T)},
\]
is continuous, i.e.,
\begin{equation} \label{ODE:atilde_stetig}
  |\widetilde{a}_\mu(u,v)| \, = \,
  | \langle \Box_\mu \widetilde{u} , \mathcal E v \rangle_{(-T,T)}|
  \leq \| \Box_\mu \widetilde{u} \|_{[H^1_0(-T,T)]'}
  \| \partial_t v \|_{L^2(0,T)}
\end{equation}
for all $u \in \mathcal H_{0,}(0,T)$ and $v \in H^1_{,0}(0,T)$, and
fulfills the inf--sup condition
\begin{equation} \label{ODE:atilde_infsup}
  \| \Box_\mu \widetilde{u} \|_{[H^1_0(-T,T)]'} =
  \sup\limits_{0 \neq v \in H^1_{,0}(0,T)}
  \frac{| \langle \Box_\mu \widetilde{u} , \mathcal E v \rangle_{(-T,T)} |}
  {\| \partial_t v \|_{L^2(0,T)}} =
  \sup\limits_{0 \neq v \in H^1_{,0}(0,T)}
  \frac{| \widetilde{a}_\mu(u,v) |}{\| \partial_t v \|_{L^2(0,T)}}
\end{equation}
for all $u \in \mathcal H_{0,}(0,T)$, where the norm representation
\eqref{ODE:BoxNormdarstellung} is used. In addition,
Lemma~\ref{ODE:Lem:Darstellung_Bilinearform} yields the representation
\begin{equation} \label{ODE:atilde_ist_a}
    \widetilde{a}_\mu(u,v) = a_\mu(u,v)
\end{equation}
for all $u \in H^1_{0,}(0,T) \subset \mathcal H_{0,}(0,T)$,
$v \in H^1_{,0}(0,T)$, which is used in the following lemma.

\begin{lemma} \label{ODE:Lem:atilde_3Bed}
  For all $0 \neq v \in H^1_{,0}(0,T)$, there exists a function
  $u_v \in \mathcal H_{0,}(0,T)$ such that
  \[
    \widetilde{a}_\mu(u_v,v) > 0 \, .
  \]
\end{lemma}

\proof{For $0 \neq v \in H^1_{,0}(0,T)$, there exists the unique solution
  $u_v \in H^1_{0,}(0,T)$ satisfying
  \[
    a_\mu(u_v,w) = \langle v , w \rangle_{L^2(0,T)} \quad
    \text{for all } w \in H^1_{,0}(0,T).
  \]
  Using the representation \eqref{ODE:atilde_ist_a}, this gives
  \[
    \widetilde{a}_\mu(u_v,w) =
    \langle v , w \rangle_{L^2(0,T)} \quad
    \text{for all } w \in H^1_{,0}(0,T),
  \]
  and in particular for $w=v$, we obtain
  \[
  \widetilde{a}_\mu(u_v,v) = \| v \|_{L^2(0,T)}^2 > 0,
\]
i.e., the assertion.}

\noindent
Next, we state the new variational setting for the scalar ordinary
differential equation \eqref{ODE}. For given $f \in [H^1_{,0}(0,T)]'$,
we consider the variational formulation to find
$u \in \mathcal H_{0,}(0,T)$ such that
\begin{equation}\label{ODE:VF_verallgemeinert}
  \widetilde{a}_\mu(u,v) = \langle f, v \rangle_{(0,T)}  \quad
  \text{for all } v \in H^1_{,0}(0,T),
\end{equation}
i.e., the operator equation
\begin{equation*}
  \mathcal E' \Box_\mu \widetilde{u} = f \quad \text{ in }
  [H^1_{,0}(0,T)]'.
\end{equation*}
With the properties of the bilinear form $\widetilde{a}_\mu(\cdot,\cdot)$,
the unique solvability of the variational formulation
\eqref{ODE:VF_verallgemeinert}, i.e., the main theorem of this section,
is proven.

\begin{theorem} \label{ODE:Thm:ExistenzEindeutigkeit}
  For each given $f \in [H^1_{,0}(0,T)]'$, there exists a unique solution
  $u \in \mathcal H_{0,}(0,T)$ of the variational formulation
  \eqref{ODE:VF_verallgemeinert}. Furthermore,
  \begin{equation*}
    \mathcal L_\mu \colon \, [H^1_{,0}(0,T)]' \to \mathcal H_{0,}(0,T),
    \qquad \mathcal L_\mu f := u,
  \end{equation*}
  is an isomorphism satisfying 
  \begin{equation*}
    \| \Box_\mu \widetilde{u} \|_{[H^1_0(-T,T)]'} =
    \| \Box_\mu \widetilde{\mathcal L_\mu f} \|_{[H^1_0(-T,T)]'} =
    \| f \|_{[H^1_{,0}(0,T)]'}.
  \end{equation*}
\end{theorem}

\proof{With the help of the Banach--Ne\v{c}as--Babu\v{s}ka theorem
  \cite[Theorem 2.6]{ErnGuermond2004}, the results in
  \eqref{ODE:atilde_stetig}, \eqref{ODE:atilde_infsup} and
  Lemma~\ref{ODE:Lem:atilde_3Bed} yield the existence and uniqueness of
  a solution $u \in \mathcal H_{0,}(0,T)$. In addition, with the
  variational formulation \eqref{ODE:VF_verallgemeinert}, the equalities
  \begin{equation*}
    \| \Box_\mu \widetilde{u} \|_{[H^1_0(-T,T)]'} =
    \sup\limits_{0 \neq v \in H^1_{,0}(0,T)}
    \frac{| \widetilde{a}_\mu(u,v) |}{\| \partial_t v \|_{L^2(0,T)}} =
    \sup\limits_{0 \neq v \in H^1_{,0}(0,T)}
    \frac{| \langle f, v \rangle_{(0,T)} |}
    {\| \partial_t v \|_{L^2(0,T)}} = \| f \|_{[H^1_{,0}(0,T)]'}
\end{equation*}
hold true, and therefore, the assertion.}

\noindent
While the unique solution $u$ of the variational formulation 
\eqref{ODE:VF_verallgemeinert} is considered in $\mathcal H_{0,}(0,T)$,
it turns out that indeed $u \in H^1_{0,}(0,T)$. In fact, the following
lemma clarifies the relation between $\mathcal H_{0,}(0,T)$ and
$H^1_{0,}(0,T)$.

\begin{lemma} \label{ODE:Lem:HNeu_ist_H10}
  There holds $\mathcal H_{0,}(0,T) = H^1_{0,}(0,T)$
  with the norm equivalence inequalities
  \begin{equation*}
    \left( 1 + \frac{4}{\pi^2} \mu T^2 \right)^{-1}
    \| \Box_\mu \widetilde{u} \|_{[H^1_0(-T,T)]'} \, \leq \, 
    \| \partial_t u \|_{L^2(0,T)} \, \leq \,
    \left( 1 + \frac{1}{2} \sqrt{\mu} T \right)
    \| \Box_\mu \widetilde{u} \|_{[H^1_0(-T,T)]'}
\end{equation*}
for all $u \in H^1_{0,}(0,T)$.
\end{lemma}

\proof{We first prove that $\mathcal H_{0,}(0,T) = H^1_{0,}(0,T)$.
  As $\mathcal H_{0,}(0,T) \supset H^1_{0,}(0,T)$, see
  Lemma~\ref{ODE:Lem:Darstellung_Bilinearform}, it remains to prove
  that $\mathcal H_{0,}(0,T) \subset H^1_{0,}(0,T)$. For this purpose,
  let $u \in \mathcal H_{0,}(0,T)$ be fixed. Consider the unique
  solution $\hat u \in H^1_{0,}(0,T)$ of the variational formulation
  \eqref{ODE:VF} for the right--hand side
  $f = \mathcal L_\mu^{-1} u \in [H^1_{,0}(0,T)]'$, where
  $\mathcal L_\mu$ is the solution operator of
  Theorem~\ref{ODE:Thm:ExistenzEindeutigkeit}. So, using
  Lemma~\ref{ODE:Lem:Darstellung_Bilinearform} and the variational
  formulations \eqref{ODE:VF}, \eqref{ODE:VF_verallgemeinert} this yields
  \begin{equation*}
    \widetilde{a}_{\mu}(\hat u,v) = a_{\mu}(\hat u,v)
    = \langle f, v \rangle_{(0,T)} = \widetilde{a}_{\mu}(u,v)
  \end{equation*}
  for all $v \in H^1_{,0}(0,T)$, i.e., $u=\hat u \in H^1_{0,}(0,T)$.
  Thus, we have $\mathcal H_{0,}(0,T) \subset H^1_{0,}(0,T)$.
  The upper norm equivalence inequality is proven by
  \begin{align*}
    \| \Box_\mu \widetilde{u} \|_{[H^1_0(-T,T)]'}
    &= \sup\limits_{0 \neq v \in H^1_{,0}(0,T)}
      \frac{| \widetilde{a}_\mu(u,v) |}{\| \partial_t v \|_{L^2(0,T)}} \\
    &= \sup\limits_{0 \neq v \in H^1_{,0}(0,T)}
      \frac{| a_\mu(u,v) |}{\| \partial_t v \|_{L^2(0,T)}} 
      \geq \frac{2}{2+\sqrt{\mu} T} \, \| \partial_t u \|_{L^2(0,T)}
\end{align*}
for all $u \in \mathcal H_{0,}(0,T) = H^1_{0,}(0,T)$, where the inf--sup
conditions \eqref{ODE:atilde_infsup}, \eqref{ODE:a_infsup} are used. The
lower inequality follows from \eqref{ODE:NormBoxuDurchAbleitungu}.}

\begin{corollary}  \label{ODE:Kor:VF_gleich}
  For all $u \in \mathcal H_{0,}(0,T)$ and all $v \in H^1_{,0}(0,T)$, the
  equality
  \begin{equation*}
    \widetilde{a}_\mu(u,v) = a_\mu(u,v)
  \end{equation*}
  is valid, i.e., the variational formulations \eqref{ODE:VF} and
  \eqref{ODE:VF_verallgemeinert} are equivalent.    
\end{corollary}

\proof{The assertion follows immediately from
  Lemma~\ref{ODE:Lem:HNeu_ist_H10} and \eqref{ODE:atilde_ist_a}.}

\begin{remark}
  Functions $u \in C^2([0,T])$ with $u(0) = 0$ are contained in
  $\mathcal H_{0,}(0,T)$, since such functions are in $H^1_{0,}(0,T)$.
  Note that the second initial condition
  $\partial_t u(t)_{|t=0} = 0$ is not incorporated in the ansatz space
  $\mathcal H_{0,}(0,T)$.
\end{remark}

\begin{remark}
  The function $u$, defined by $u(t) = \sin (\sqrt{\mu} t)$ for
  $t \in (0,T),$ is obviously in $H^1_{0,}(0,T)$ and so, in
  $\mathcal H_{0,}(0,T)$. For this function, we have
  \[
    \| \partial_{tt} u + \mu u \|_{[H^1_{,0}(0,T)]'}=0 .
  \]
  For
  \[
    \widetilde{u}(t) =
    \begin{cases}
     0 & \text{for } t \in (-T, 0], \\
     \sin (\sqrt{\mu} t) & \text{for } t \in (0,T),
     \end{cases}
   \]
   the first--order distributional derivative is identified with the function
   \[
     \partial_t \widetilde{u}(t) =
     \begin{cases}
     0 & \text{for } t \in (-T, 0], \\
     \sqrt{\mu} \cos (\sqrt{\mu} t) & \text{for } t \in (0,T),
     \end{cases}
   \]
   i.e., it is a regular distribution. To compute the second--order
   distributional derivative of $\widetilde{u}$, we consider
   \begin{align*}
     \langle \partial_{tt} \widetilde{u}, \varphi \rangle_{(-T,T)}
     &= -\int_{-T}^T \partial_t \widetilde{u}(t) \, \partial_t \varphi(t)
       \, \mathrm dt \\
     &= -\int_0^T \sqrt{\mu} \cos (\sqrt{\mu} t) \,
       \partial_t \varphi(t) \, \mathrm dt \\ 
     &= -\sqrt{\mu} \cos (\sqrt{\mu} t) \, \varphi(t) \Big|_0^T +
       \int_0^T \Big( - \mu \sin (\sqrt{\mu} t) \Big) \, \varphi(t)
       \, \mathrm dt \\
     &=  \sqrt{\mu} \, \varphi(0) -
       \mu \langle \widetilde{u}, \varphi \rangle_{(-T,T)}
    \end{align*}
    for all $\varphi \in C^\infty_0(-T,T).$ Hence,
    \[
      \Box_\mu \widetilde{u} =
      \partial_{tt} \widetilde{u} + \mu \widetilde{u} =
      \sqrt{\mu} \, \delta_0
    \]
    is a singular distribution with the Dirac distribution
    $\delta_0 \in H^{-1}_{[0,T]}(-T,T) \subset [H^1_0(-T,T)]'$.
    Furthermore, it follows that
    \begin{align*}
      \| \Box_\mu \widetilde{u} \|_{[H^1_0(-T,T)]'}
      &= \sup\limits_{0 \neq v \in H^1_{,0}(0,T)}
        \frac{| \langle \Box_\mu \widetilde{u},\mathcal E v \rangle_{(-T,T)} |}
        {\| \partial_t v \|_{L^2(0,T)}} \, = \,
        \sup\limits_{0 \neq v \in H^1_{,0}(0,T)}
        \frac{\sqrt{\mu} \, |\langle \delta_0,\mathcal E v  \rangle_{(-T,T)}| }
        {\| \partial_t v \|_{L^2(0,T)}} \\
      &= \sup\limits_{0 \neq v \in H^1_{,0}(0,T)}
        \frac{\sqrt{\mu} \, |v(0)| }{\| \partial_t v \|_{L^2(0,T)}} > 0
    \end{align*}
    for, e.g., $v(t) = \frac{1}{T}(T-t)$, where the norm representation
    \eqref{ODE:BoxNormdarstellung} is used. To summarize, the function
    $u(t) = \sin (\sqrt{\mu} t)$ for $t \in (0,T)$ with
    \begin{equation*}
      \| \Box_\mu \widetilde{u} \|_{[H^1_0(-T,T)]'} > 0 \quad
      \text{ and } \quad \| \partial_{tt} u + \mu u \|_{[H^1_{,0}(0,T)]'}=0
    \end{equation*}
    solves the variational formulations \eqref{ODE:VF} and
    \eqref{ODE:VF_verallgemeinert} for the right--hand side
    $f_{v_0} \in [H^1_{,0}(0,T)]'$,
    \begin{equation*}
      \langle f_{v_0}, v \rangle_{(0,T)} =
      \sqrt{\mu} v(0), \quad v \in H^1_{,0}(0,T),
    \end{equation*}
    which realizes the initial condition
    $\partial_t u(t)_{|t=0} = v_0 := \sqrt{\mu}$.
\end{remark}

\begin{remark}
    The variational formulation \eqref{ODE:VF_verallgemeinert} is the weak
    formulation of the differential equation
    \[
        \partial_{tt} \widetilde{u}(t) + \mu \widetilde{u}(t)
        = \widetilde{f}(t)= \begin{cases}
            0    & \text{for } t \in (-T,0),               \\
            f(t) & \text{for } t \in [0,T),
          \end{cases}
    \]
    which can be written as coupled system, using $u(t)=\widetilde{u}(t)$
    for $t \in (0,T)$, and $u_-(t)=\widetilde{u}(t)$ for $ t \in (-T,0)$,
    \begin{align*}
      \partial_{tt} u(t) + \mu u(t)
      &= f(t) \quad  &&\text{for } t \in (0,T), \\
      \partial_{tt} u_-(t) + \mu u_-(t)
      &= 0 \quad  &&\text{for } t \in (-T,0), \quad
                     u_-(-T)= \partial_t u_-(t)_{|t=-T} =0
    \end{align*}
    together with the transmission interface conditions
    \[
      u(0) = u_-(0), \quad
      \partial_t u(t)_{|t=0} -  \partial_t u_-(t)_{|t=0} = v_0
    \]
    with given $v_0 \in \R$, satisfying
    $\langle \widetilde{f}, z \rangle_{(-T,T)} = v_0 z(0)$ for all
    $z \in H^1_0(-T,T).$ The conditions
    $u_-(-T)= \partial_t u_-(t)_{|t=-T} =0$ lead to
    $u_-(t) = \widetilde{u}(t) = 0 $ for $t \in (-T,0)$, which
    finally implies the initial conditions
    \[
        u(0) = 0, \quad \partial_t u(t)_{|t=0} = v_0.
    \]
\end{remark}

\section{A generalized variational formulation for the wave
  equation} \label{Sec:Welle}
In this section, we generalize the approach, as introduced for the solution
of the ordinary differential equation \eqref{ODE}, to end up with a generalized
inf--sup stable variational formulation for the Dirichlet boundary value
problem for the wave equation \eqref{Welle}. For this purpose, we first
introduce notations analogously to them of Section~\ref{Sec:ODE}.

In addition to the space--time domain $Q:=\Omega \times (0,T)$ we consider
the extended domain $Q_-:= \Omega \times (-T,T)$. The dual space
$[H^{1,1}_{0;\,,0}(Q)]'$ is characterized as completion of $L^2(Q)$ with
respect to the Hilbertian norm
\begin{equation*}
  \| f \|_{[H^{1,1}_{0;\,,0}(Q)]'} =
  \sup_{0 \neq v \in H^{1,1}_{0;\,,0}(Q) } 
  \frac{| \langle f, v \rangle_Q |}{ \| v \|_{H^{1,1}_{0;\,,0}(Q)}},
\end{equation*}
where $\langle \cdot , \cdot \rangle_Q$ denotes the duality pairing as
extension of the inner product in $L^2(Q)$. Note that
$[H^{1,1}_{0;\,,0}(Q)]'$ is a Hilbert space, see Section~\ref{Sec:ODE}.
For given $u \in L^2(Q)$, we define the extension
$\widetilde{u} \in L^2(Q_-)$ by
\begin{equation}\label{Welle:DefNullFortsetzung}
  \widetilde{u}(x,t) :=
  \begin{cases}
    u(x,t) & \text{for } (x,t) \in Q, \\
    0 & \text{for } (x,t) \in Q_- \setminus Q.
  \end{cases}
\end{equation}
The application of the wave operator $\Box := \partial_{tt} - \Delta_x$
to $\widetilde{u}$ is defined as a distribution on $Q_-$, i.e., for 
all test functions $\varphi \in C^\infty_0(Q_-)$, we define
\begin{equation}\label{Welle:Def_verallgemeinerte_Ableitung}
  \langle \Box \widetilde{u} , \varphi \rangle_{Q_-} :=
  \int_{-T}^T \int_\Omega \widetilde{u}(x,t) \, \Box \varphi(x,t)
  \, \mathrm dx \, \mathrm dt = \int_0^T \int_\Omega u(x,t) \,
  \Box \varphi(x,t) \, \mathrm dx \, \mathrm dt.
\end{equation}
This motivates to consider the dual space $[H^1_0(Q_-)]'$ of $H^1_0(Q_-)$,
which is characterized as completion of $L^2(Q_-)$ with respect to the
Hilbertian norm
\begin{equation*}
  \| g \|_{[H^1_0(Q_-)]'} := \sup_{0 \neq z \in H^1_0(Q_-) } 
  \frac{|\langle g, z \rangle_{Q_-}|}{ \| z \|_{H^1_0(Q_-)}},
\end{equation*}
where the inner product
\begin{equation*}
  \langle z_1, z_2 \rangle_{H^1_0(Q_-)} =
  \langle \partial_t z_1, \partial_t z_2 \rangle_{L^2(Q_-)} +
  \langle \nabla_x z_1, \nabla_x z_2 \rangle_{L^2(Q_-)},
  \quad z_1, z_2 \in H^1_0(Q_-),
\end{equation*}
induces the norm $\| \cdot \|_{H^1_0(Q_-)}$, and
$\langle \cdot , \cdot \rangle_{Q_-}$ denotes the duality pairing as
extension of the inner product in $L^2(Q_-)$,
see \cite[Satz 17.3]{Wloka1982}. Note that $[H^1_0(Q_-)]'$ is a Hilbert
space, see Section~\ref{Sec:ODE}. In addition we define the subspace
\begin{equation*}
  H^{-1}_{|\overline{Q}}(Q_-) := \Big \{ g \in [H^1_0(Q_-)]' \colon
  \forall z \in H^1_0(Q_-) \text{ with }
  \mathrm{supp}\, z \subset \Omega \times (-T,0) \colon  \,
  \langle g, z \rangle_{Q_-} = 0 \Big\}
\end{equation*}
of $[H^1_0(Q_-)]'$, endowed with the Hilbertian norm
$\norm{\cdot}_{[H^1_0(Q_-)]'}$. To characterize the
subspace $H^{-1}_{|\overline{Q}}(Q_-)$, we introduce the following notations.
Let ${\mathcal{R}} \colon H^1_0(Q_-) \to H^{1,1}_{0;,0}(Q)$ be the continuous and
surjective restriction operator, defined by $\mathcal R z = z_{|Q}$ for
$z \in H^1_0(Q_-)$, with its adjoint operator
$\mathcal R' \colon [H^{1,1}_{0;,0}(Q)]' \to [H^1_0(Q_-)]'$. Furthermore, let
$\mathcal E \colon H^{1,1}_{0;,0}(Q) \to H^1_0(Q_-)$ be any continuous and
injective extension operator with its adjoint operator
$\mathcal E' \colon [H^1_0(Q_-)]' \to [H^{1,1}_{0;,0}(Q)]'$, satisfying
\begin{equation*}
  \| \mathcal E v \|_{H^1_0(Q_-)} \leq c_{\mathcal E} \, \| v \|_{H^{1,1}_{0;,0}(Q)}
\end{equation*}
with a constant $c_{\mathcal E} > 0$, and $\mathcal R \mathcal E v = v$ for
all $v \in H^{1,1}_{0;,0}(Q)$. An example for such an extension operator is
given by reflection in $\Omega \times \{0\}$, i.e., consider the
function $\overline{v}$, defined by
\begin{equation*}
  \overline{v}(x,t) =
  \begin{cases}
    v(x,t)  & \text{for } (x,t) \in \Omega \times [0,T), \\
    v(x,-t) & \text{for } (x,t) \in \Omega \times (-T,0)
  \end{cases}
\end{equation*}
for $(x,t) \in Q_-$, and a given function $v \in H^{1,1}_{0;,0}(Q)$, which
leads to a constant $c_{\mathcal E} = 2$ in this particular case.
With this, we prove the following lemma as the counter part of
Lemma \ref{ODE:Lem:DualH10Support}.

\begin{lemma} \label{Welle:Lem:DualH10Support}
  The spaces $(H^{-1}_{|\overline{Q}}(Q_-), \| \cdot \|_{[H^1_0(Q_-)]'})$ and
  $([H^{1,1}_{0;,0}(Q)]', \| \cdot \|_{[H^{1,1}_{0;,0}(Q)]'})$ are isometric,
  i.e., the mapping
  \begin{equation*}
    \mathcal E'_{|H^{-1}_{|\overline{Q}}(Q_-)} \colon \,
    H^{-1}_{|\overline{Q}}(Q_-) \to  [H^{1,1}_{0;,0}(Q)]'
  \end{equation*}
  is bijective with
  \begin{equation*}
    \| g \|_{[H^1_0(Q_-)]'} = \| \mathcal E' g \|_{[H^{1,1}_{0;,0}(Q)]'}
    \quad \text{ for all } g \in H^{-1}_{|\overline{Q}}(Q_-).
  \end{equation*}
  In addition, for $g \in H^{-1}_{|\overline{Q}}(Q_-)$, the relation
  \begin{equation} \label{Welle:Lem:DualH10Support:gDurchR}
    \langle g, z \rangle_{Q_-} =
    \langle \mathcal E' g, \mathcal R z \rangle_Q
    \quad \mbox{for all} \; z \in H^1_0(Q_-)
    \end{equation}
    i.e., $\mathcal R' \mathcal E' g = g$, holds true. In particular,
    the subspace $H^{-1}_{|\overline{Q}}(Q_-) \subset [H^1_0(Q_-)]'$ is
    closed, i.e., complete.
\end{lemma}

\proof{First, we prove that
  $\| g \|_{[H^1_0(Q_-)]'} = \| \mathcal E' g \|_{[H^{1,1}_{0;,0}(Q)]'}$ for
  all functionals $g \in H^{-1}_{|\overline{Q}}(Q_-)$. For this purpose,
  let $g \in H^{-1}_{|\overline{Q}}(Q_-)$ be arbitrary but fixed. The Riesz
  representation theorem gives the unique element $z_g \in H^1_0(Q_-)$ with
  \begin{equation*}
    \langle g, z \rangle_{Q_-} = \langle z_g, z \rangle_{H^1_0(Q_-)}
    \quad \text{ for all } z \in H^1_0(Q_-),
  \end{equation*}
  and $\| g \|_{[H^1_0(Q_-)]'} = \| z_g \|_{H^1_0(Q_-)}$. It holds true
  that $z_{g|\Omega \times (-T,0)} = 0$, since we have
  \begin{equation*}
    0 = \langle g, z \rangle_{Q_-} =
    \langle  z_g,  z \rangle_{H^1_0(Q_-)} =
    \int_{-T}^0 \int_\Omega \Big[
    \partial_t z_g(x,t) \, \partial_t z(x,t) +
    \nabla_x z_g(x,t) \cdot \nabla_x z(x,t) \Big] \, \mathrm dx \, \mathrm dt
  \end{equation*}
  for all $z \in H^1_0(Q_-)$ with
  $\mathrm{supp}\, z \subset \Omega \times (-T,0)$. Hence, we have
  \begin{equation} \label{Welle:Lem:DualH10Support:Beweis:gDurchR}
    \langle g, z \rangle_{Q_-} = \langle z_g, z \rangle_{H^1_0(Q_-)} =
    \langle \mathcal R z_g, \mathcal R z \rangle_{H^{1,1}_{0;,0}(Q)}
  \end{equation}
  for all $z \in H^1_0(Q_-)$. So, using
  \eqref{Welle:Lem:DualH10Support:Beweis:gDurchR} with $z=\mathcal E v$
  for $v \in H^{1,1}_{0;,0}(Q)$ this gives
  \begin{equation} \label{Welle:Lem:DualH10Support:Beweis:E'gDurchR}
    \langle \mathcal E'g, v \rangle_Q =
    \langle g, \mathcal E v \rangle_{Q_-} =
    \langle \mathcal R z_g, \mathcal R \mathcal E v \rangle_{H^{1,1}_{0;,0}(Q)}
    = \langle \mathcal R z_g, v \rangle_{H^{1,1}_{0;,0}(Q)},
  \end{equation}
  i.e.,
  \begin{equation*}
    \| \mathcal E'g \|_{[H^{1,1}_{0;,0}(Q)]'} =
    \| \mathcal R z_g \|_{H^{1,1}_{0;,0}(Q)} =
    \| z_g \|_{H^1_0(Q_-)} =
    \| g \|_{[H^1_0(Q_-)]'} .
  \end{equation*}
  Next, we prove that $\mathcal E'_{|H^{-1}_{|\overline{Q}}(Q_-)} $ is
  surjective. For this purpose, let $f \in [H^{1,1}_{0;,0}(Q)]'$ be given.
  Set $g_f = \mathcal R'f$, i.e.,
  \begin{equation*}
    \langle g_f, z \rangle_{Q_-} =
    \langle \mathcal R' f , z \rangle_{Q_-} =
    \langle f, \mathcal R z \rangle_Q 
  \end{equation*}
  for all $z \in H^1_0(Q_-)$. With this it follows immediately that
  $g_f \in H^{-1}_{|\overline{Q}}(Q_-)$. Moreover, we have
  \begin{equation*}
    \langle \mathcal E' g_f, v \rangle_Q =
    \langle  g_f, \mathcal E v \rangle_{Q_-} =
    \langle f, \mathcal R \mathcal E v \rangle_Q =
    \langle f, v \rangle_Q
\end{equation*}
for all $v \in H^{1,1}_{0;,0}(Q)$, i.e., $\mathcal E' g_f = f$ in 
$[H^{1,1}_{0;,0}(Q)]'$. In other words,
$\mathcal E'_{|H^{-1}_{|\overline{Q}}(Q_-)} $ is surjective.
Finally, \eqref{Welle:Lem:DualH10Support:gDurchR} follows from
\eqref{Welle:Lem:DualH10Support:Beweis:gDurchR} and
\eqref{Welle:Lem:DualH10Support:Beweis:E'gDurchR} for $v = \mathcal R z$
for any $z \in H^1_0(Q_-)$. The last assertion of the lemma is
straightforward.}

\noindent
The last lemma gives immediately the following corollary.

\begin{corollary} \label{Welle:Kor:NormHSupp0T}
  For all $g \in H^{-1}_{|\overline{Q}}(Q_-)$, the norm representation
  \begin{equation*} 
    \| g \|_{[H^1_0(Q_-)]'} =
    \sup\limits_{0 \neq v \in H^{1,1}_{0;,0}(Q)}
    \frac{|\langle g, \mathcal E v \rangle_{Q_-}|}{\| v \|_{H^{1,1}_{0;,0}(Q)}}
  \end{equation*}
  holds true.
\end{corollary}

\proof{Let $g \in H^{-1}_{|\overline{Q}}(Q_-)$ be arbitrary but fixed.
  With Lemma~\ref{Welle:Lem:DualH10Support}, we have
  \begin{equation*}
    \| g \|_{[H^1_0(Q_-)]'} =
    \| \mathcal E' g \|_{[H^{1,1}_{0;,0}(Q)]'} =
    \sup\limits_{0 \neq v \in H^{1,1}_{0;,0}(Q)}
    \frac{|\langle \mathcal E' g , v\rangle_Q|}{\| v \|_{H^{1,1}_{0;,0}(Q)}} =
    \sup\limits_{0 \neq v \in H^{1,1}_{0;,0}(Q)}
    \frac{|\langle g , \mathcal E v\rangle_{Q_-}|}{\| v \|_{H^{1,1}_{0;,0}(Q)}},
\end{equation*}
i.e., the assertion is proven.}

\noindent
Next, we introduce
\[
  \mathcal H (Q) := \Big \{ u= \widetilde{u}_{|Q} : 
  \widetilde{u} \in L^2(Q_-), \;
  \widetilde{u}_{|\Omega \times (-T,0)} = 0, \;
  \Box \widetilde{u} \in [H^1_0(Q_-)]' \Big \},
\]
with the norm
\[
  \| u \|_{\mathcal{H}(Q)} :=
  \sqrt{\| u \|^2_{L^2(Q)} + \| \Box \widetilde{u} \|^2_{[H^1_0(Q_-)]'}} \; .
\]
For a function $u \in \mathcal{H}(Q)$, the condition
$\Box \widetilde{u} \in [H^1_0(Q_-)]'$ involves that there exists an
element $f_u \in [H^1_0(Q_-)]'$ with
\begin{equation*}
  \langle \Box \widetilde{u} , \varphi \rangle_{Q_-}  =
  \spf{f_u}{\varphi}_{Q_-} \quad \text{ for all } \varphi \in  C^\infty_0(Q_-).
\end{equation*}
Note that $\varphi \in H^1_0(Q_-)$ for $\varphi \in C^\infty_0(Q_-)$, and
that $C^\infty_0(Q_-)$ is dense in $H^1_0(Q_-)$. Hence, the element
$f_u \in [H^1_0(Q_-)]'$ is unique and therefore, in the following, we
identify the distribution
$\Box \widetilde{u} \colon \, C^\infty_0(Q_-) \to \R$ with the
functional $f_u \colon \, H^1_0(Q_-) \to \R$.

Next, we state properties of the space $\mathcal H (Q)$. Clearly,
$(\mathcal H (Q),\| \cdot \|_{\mathcal H (Q)})$ is a normed
vector space and it is even a Banach space.

\begin{lemma}\label{Welle:Lem:Banachraum}
  The normed vector space $(\mathcal H (Q),\| \cdot \|_{\mathcal H (Q)})$
  is a Banach space.
\end{lemma}

\proof{Consider a Cauchy sequence
$(u_n)_{n \in {\N }} \subset \mathcal H (Q)$. Hence,
$(u_n)_{n \in {\N }} \subset L^2(Q)$ is also a Cauchy sequence in $L^2(Q)$,
and $(\Box \widetilde{u}_n)_{n \in {\N }} \subset [H^1_0(Q_-)]'$ is also
a Cauchy sequence in $[H^1_0(Q_-)]'$. So, there exist $u \in L^2(Q)$
and $ f \in [H^1_0(Q_-)]'$ with
\[
  \lim\limits_{n \to \infty} \| u_n - u \|_{L^2(Q)} = 0, \quad
  \lim\limits_{n \to \infty} \| \Box \widetilde{u}_n - f \|_{[H^1_0(Q_-)]'} = 0.
\]
For $\varphi \in C_0^\infty(Q_-)$, we have
\begin{align*}
  \langle \Box \widetilde{u}, \varphi \rangle_{Q_-}
  &=
    \langle \widetilde{u} , \Box \varphi \rangle_{L^2(Q_-)} =
    \int_0^T \int_\Omega \, u(x,t) \Box \varphi(x,t) \,
    \mathrm dx \, \mathrm dt \\
  &= \lim\limits_{n \to \infty} \int_0^T \int_\Omega u_n(x,t) \,
    \Box \varphi(x,t) \, \mathrm dx \, \mathrm dt \\ 
  &= \lim\limits_{n \to \infty} \langle \widetilde{u}_n ,
    \Box \varphi \rangle_{L^2(Q_-)} =
    \lim\limits_{n \to \infty} \langle \Box \widetilde{u}_n ,
    \varphi \rangle_{Q_-} \, = \, \langle f , \varphi \rangle_{Q_-},
\end{align*}
i.e., $\Box \widetilde{u} = f \in [H^1_0(Q_-)]'$. Hence,
$u \in \mathcal H (Q)$ follows.}

\noindent
With the abstract inner product $\langle \cdot, \cdot \rangle_{[H^1_0(Q_-)]'}$
of $[H^1_0(Q_-)]'$, the inner product
\begin{equation*}
  \langle \cdot, \cdot \rangle_{\mathcal H (Q)} :=
  \langle \cdot, \cdot \rangle_{L^2 (Q)} +
  \langle \Box \widetilde{(\cdot)},
  \Box \widetilde{(\cdot)} \rangle_{[H^1_0(Q_-)]'}
\end{equation*}
induces the norm $\| \cdot \|_{\mathcal H (Q)}$. Hence, the space
$(\mathcal H (Q), \langle \cdot, \cdot \rangle_{\mathcal H (Q)})$ is even
a Hilbert space, but this abstract inner product is not used explicitly
in the remainder of this work.

\begin{lemma} \label{Welle:Lem:BoxInH10Support}
  For all $u \in \mathcal H (Q)$ there holds
  $\Box \widetilde{u} \in H^{-1}_{|\overline{Q}}(Q_-)$ with
  \begin{equation} \label{Welle:BoxNormdarstellung}
    \| \Box \widetilde{u} \|_{[H^1_0(Q_-)]'} =
    \sup\limits_{0 \neq v \in H^{1,1}_{0;,0}(Q)}
    \frac{|\langle \Box \widetilde{u}, \mathcal E v \rangle_{Q_-}|}
    {\| v \|_{H^{1,1}_{0;,0}(Q)}} .
    \end{equation}
  \end{lemma}
  
\proof{First, we prove that $\Box \widetilde{u} \in H^{-1}_{|\overline{Q}}(Q_-)$.
  For this purpose, let $u \in \mathcal H (Q)$ and $z \in H^1_0(Q_-)$ with
  $\mathrm{supp}\, z \subset \Omega \times (-T,0)$ be arbitrary but fixed.
  Due to $z_{|\Omega \times (-T,0)} \in H^1_0(\Omega \times (-T,0))$ there
  exists a sequence
  $(\psi_n)_{n \in \N} \subset C^\infty_0(\Omega \times (-T,0))$ with
  $\| z_{|\Omega \times (-T,0)} - \psi_n \|_{H^1_0(\Omega \times (-T,0))} \to 0$
  as $n \to \infty$, where
  \begin{equation*}
    \| w \|_{H^1_0(\Omega \times (-T,0))} =
    \left( \int_{-T}^0 \int_\Omega \Big[
    |\partial_t w(x,t)|^2 + |\nabla_x w(x,t)|^2 \Big] \,
    \mathrm dx \, \mathrm dt \right)^{1/2}
    \end{equation*}
    for $w \in H^1_0(\Omega \times (-T,0))$. For $n \in \N$, define
    \begin{equation*}
      \varphi_n(x,t) =
      \begin{cases}
        \psi_n(x,t) & \text{for } (x,t) \in \Omega \times (-T,0), \\
        0       & \text{for } (x,t) \in \Omega \times [0,T),
      \end{cases}
    \end{equation*}
    i.e., $(\varphi_n)_{n \in \N} \subset C^\infty_0(Q_-)$ satisfies
    \begin{equation*}
      \|  z - \varphi_n \|_{H^1_0(Q_-)} =
      \| z_{|\Omega \times (-T,0)} - \psi_n \|_{H^1_0(\Omega \times (-T,0))} \to 0
    \end{equation*}
    as $n \to \infty$. So, it follows that
    \begin{equation*}
      \langle \Box \widetilde{u}, z \rangle_{Q_-} =
      \lim_{n\to\infty} \langle \Box \widetilde{u}, \varphi_n \rangle_{Q_-} =
      \lim_{n \to \infty} \int_0^T \int_\Omega u(x,t) \,
      \Box \varphi_n(x,t) \, \mathrm dx \, \mathrm dt = 0
    \end{equation*}
    and therefore, the assertion. The norm representation follows from
    $\Box \widetilde{u} \in H^{-1}_{|\overline{Q}}(Q_-)$ and
    Corollary~\ref{Welle:Kor:NormHSupp0T}.}

\begin{lemma}\label{Welle:Lem:Darstellung_Bilinearform}
  It holds true that $H^{1,1}_{0;0,}(Q) \subset \mathcal H (Q)$. Furthermore,
  each $u \in H^{1,1}_{0;0,}(Q)$ with zero extension $\widetilde{u}$, as
  defined in \eqref{Welle:DefNullFortsetzung}, satisfies
  \begin{equation} \label{Welle:NormBoxuDurchAbleitungu}
    \| \Box \widetilde{u} \|_{[H^1_0(Q_-)]'} \leq  \| u \|_{H^{1,1}_{0;0,}(Q)},
  \end{equation}
  and
  \begin{equation}\label{Welle:Darstellung_Bilinearform}
    \langle \Box \widetilde{u} , z \rangle_{Q_-} = a(u,\mathcal R z) =
    - \langle \partial_t u , \partial_t \mathcal R z \rangle_{L^2(Q)} +
    \langle \nabla_x u , \nabla_x \mathcal R z \rangle_{L^2(Q)}
  \end{equation}
  for all $z \in H^1_0(Q_-),$ where $a(\cdot,\cdot)$ is the bilinear
  form \eqref{Welle:a}.
\end{lemma}

\proof{First, we prove that $H^{1,1}_{0;0,}(Q) \subset \mathcal H (Q)$.
  For $u \in H^{1,1}_{0;0,}(Q)$, we define the extension $\widetilde{u}$,
  see \eqref{Welle:DefNullFortsetzung}. By construction, we have
  $\widetilde{u} \in L^2(Q_-)$, and $\widetilde{u}_{|\Omega \times (-T,0)}=0$.
  It remains to prove that $\Box \widetilde{u} \in [H^1_0(Q_-)]'$. For this
  purpose, define the functional $f_u \in [H^1_0(Q_-)]'$ by
  \begin{equation*}
    \langle f_u, z \rangle_{Q_-} =  a(u,\mathcal R z)
  \end{equation*}
  for all $z \in H^1_0(Q_-)$, where $a(\cdot,\cdot)$ is the bilinear
  form \eqref{Welle:a}. The continuity of $f_u$ follows from
  \begin{equation*}
    |\langle f_u, z \rangle_{Q_-}| = |a(u,\mathcal R z)| \leq
    \| u \|_{H^{1,1}_{0;0,}(Q)} \| \mathcal R z \|_{H^{1,1}_{0;,0}(Q)} \leq
    \| u \|_{H^{1,1}_{0;0,}(Q)} \| z \|_{H^1_0(Q_-)}  
  \end{equation*}
  for all $z \in H^1_0(Q_-)$, where the estimate \eqref{Welle:a:stetig} is
  used. Using the definition (\ref{Welle:Def_verallgemeinerte_Ableitung})
  and integration by parts, this gives
  \begin{align*}
    \langle \Box \widetilde{u} , \varphi \rangle_{Q_-}
    &= \int_0^T \int_\Omega u(x,t) \, \Box \varphi(x,t) \,
      \mathrm dx \, \mathrm dt \\
    &= - \langle \partial_t u , \partial_t \mathcal R \varphi \rangle_{L^2(Q)}
      + \langle \nabla_x u , \nabla_x \mathcal R \varphi \rangle_{L^2(Q)} =
      \langle f_u, \varphi \rangle_{Q_-} 
\end{align*}
for all $\varphi \in C^\infty_0(Q_-)$, i.e.,
$\Box \widetilde{u} = f_u \in [H^1_0(Q_-)]'$. The equality
\eqref{Welle:Darstellung_Bilinearform} follows from the density
of $C^\infty_0(Q_-)$ in $H^1_0(Q_-)$. The estimate
\eqref{Welle:NormBoxuDurchAbleitungu} is proven by
\begin{align*}
  \| \Box \widetilde{u} \|_{[H^1_0(Q_-)]'}
  &= \sup\limits_{0 \neq v \in H^{1,1}_{0;,0}(Q)}
    \frac{| \langle \Box \widetilde{u}, \mathcal E v \rangle_{Q_-}|}
    {\| v \|_{H^{1,1}_{0;,0}(Q)}} =
    \sup\limits_{0 \neq v \in H^{1,1}_{0;,0}(Q)}
    \frac{| \langle f_u, \mathcal E v \rangle_{Q_-}|}
    {\| v \|_{H^{1,1}_{0;,0}(Q)}} \\
  &=\sup\limits_{0 \neq v \in H^{1,1}_{0;,0}(Q)}
    \frac{ |a(u,\mathcal R \mathcal E v) | }
    {\| v \|_{H^{1,1}_{0;,0}(Q)}} \leq \| u \|_{H^{1,1}_{0;0,}(Q)} 
\end{align*}
when using the norm representation \eqref{Welle:BoxNormdarstellung}, the
equality \eqref{Welle:Darstellung_Bilinearform}, and \eqref{Welle:a:stetig}.}

\noindent
Next, by completion, we define the Hilbert space
\[
  \mathcal H_{0,}(Q) :=
  \overline{H^{1,1}_{0;0,}(Q)}^{\| \cdot \|_{\mathcal H (Q)}}
  \subset \mathcal H (Q),
\]
endowed with the Hilbertian norm $\norm{\cdot}_{\mathcal H(Q)}$, i.e., 
\begin{equation*}
  \mathcal H_{0,}(Q) = \Big \{ v \in \mathcal H(Q) \colon \,
  \exists (v_n)_{n \in \N} \subset H^{1,1}_{0;0,}(Q) \text{ with }
  \norm{v_n - v}_{\mathcal H(Q)} \to 0 \Big \}.
\end{equation*}

\begin{lemma}\label{Welle:Lem:Friedrichs}
  For $u \in \mathcal H_{0,}(Q)$ there holds
  \[
    \| \Box \widetilde{u} \|_{[H^1_0(Q_-)]'} \geq
    \frac{\sqrt{2}}{T} \, \| u \|_{L^2(Q)} .
  \]
\end{lemma}

\proof{For $0 \neq u \in \mathcal H_{0,}(Q)$, there exists a
non--trivial sequence $(u_n)_{n \in \N } \subset H^{1,1}_{0;0,}(Q)$,
$u_n \not\equiv 0$, with
\[
  \lim\limits_{n \to \infty} \| u - u_n \|_{\mathcal H (Q)} = 0.
\]
For each $u_n \in H^{1,1}_{0;0,}(Q)$, we define $w_n \in H^{1,1}_{0;,0}(Q)$
as unique solution of the variational formulation
\[
  a(v,w_n) = \langle u_n , v \rangle_{L^2(Q)}
  \quad \text{for all } v \in H^{1,1}_{0;0,}(Q)
\]
with the bilinear form \eqref{Welle:a}. In particular for $v=u_n$, this gives
\[
  a(u_n,w_n) = \| u_n \|^2_{L^2(Q)} .
\]
Analogously to the estimate \eqref{Welle:H1Stabilitaet} for the solution of
\eqref{Welle}, we conclude
\[
  \| w_n \|_{H^{1,1}_{0;,0}(Q)} \leq \frac{1}{ \sqrt{2}} \, T \,
  \| u_n \|_{L^2(Q)} \, .
\]
For the zero extension $\widetilde{u}_n \in L^2(Q_-)$ of
$u_n \in H^{1,1}_{0;0,}(Q)$, we obtain, when using the norm representation
\eqref{Welle:BoxNormdarstellung} and \eqref{Welle:Darstellung_Bilinearform},
that
\begin{align*}
  \| \Box \widetilde{u}_n \|_{[H^1_0(Q_-)]'}
  &= \sup\limits_{0 \neq v \in H^{1,1}_{0;,0}(Q)}
    \frac{|\langle \Box \widetilde{u}_n, \mathcal E v \rangle_{Q_-}|}
    {\| v \|_{H^{1,1}_{0;,0}(Q)}} \geq
    \frac{|\langle \Box \widetilde{u}_n , \mathcal E w_n \rangle_{Q_-}|}
    {\| w_n \|_{H^{1,1}_{0;,0}(Q)}} \\
  &= \frac{| a(u_n,w_n) |}{\| w_n \|_{H^{1,1}_{0;,0}(Q)}} \, = \,
    \frac{\| u_n \|^2_{L^2(Q)}}{\| w_n \|_{H^{1,1}_{0;,0}(Q)}} \, \geq \,
    \frac{\sqrt{2}}{T} \, \| u_n \|_{L^2(Q)} ,
\end{align*}
and the assertion follows by completion for $n \to \infty$.}

\begin{corollary} 
  The inner product space $\left(\mathcal H_{0,}(Q),
    \langle \Box \widetilde{(\cdot)}, \Box \widetilde{(\cdot)}
    \rangle_{[H^1_0(Q_-)]'} \right)$ is complete, i.e., a Hilbert space.
\end{corollary}

\proof{The assertion follows immediately from
  Lemma~\ref{Welle:Lem:Friedrichs}.}

\noindent
In the following, $\mathcal H_{0,}(Q)$ is endowed with the Hilbertian norm
$\| \Box \widetilde{(\cdot)} \|_{[\overline{H}^1_0(Q_-)]'}$. With this new
Hilbert space, the bilinear form
\[
  \widetilde{a}(\cdot,\cdot) \colon \,
  \mathcal H_{0,}(Q) \times H^{1,1}_{0;,0}(Q) \to \R, \quad 
  \widetilde{a}(u,v) :=
  \langle \Box \widetilde{u} , \mathcal E v \rangle_{Q_-},
\]
is continuous, i.e.,
\begin{equation} \label{Welle:atilde_stetig}
  |\widetilde{a}(u,v)| \, = \,
  | \langle \Box \widetilde{u} , \mathcal E v \rangle_{Q_-}| \, \leq \,
  \| \Box \widetilde{u} \|_{[H^1_0(Q_-)]'} \| v \|_{H^{1,1}_{0;,0}(Q)}
\end{equation}
for all $u \in \mathcal H_{0,}(Q)$ and $v \in H^{1,1}_{0;,0}(Q)$, and fulfills
the inf--sup condition
\begin{equation} \label{Welle:atilde_infsup}
  \| \Box \widetilde{u} \|_{[H^1_0(Q_-)]'} =
  \sup\limits_{0 \neq v \in H^{1,1}_{0;,0}(Q)}
  \frac{|\langle \Box \widetilde{u} , \mathcal E v \rangle_{Q_-}|}
  {\| v \|_{H^{1,1}_{0;,0}(Q)}} =
  \sup\limits_{0 \neq v \in H^{1,1}_{0;,0}(Q)}
  \frac{| \widetilde{a}(u,v) |}{\| v \|_{H^{1,1}_{0;,0}(Q)}}
\end{equation}
for all $u \in \mathcal H_{0,}(Q)$, where the norm representation
\eqref{Welle:BoxNormdarstellung} is used. In addition,
Lemma~\ref{Welle:Lem:Darstellung_Bilinearform} yields the representation
\begin{equation} \label{Welle:atilde_ist_a}
  \widetilde{a}(u,v) = a(u,v)
\end{equation}
for all $u \in H^{1,1}_{0;0,}(Q) \subset \mathcal H_{0,}(Q)$,
$v \in H^{1,1}_{0;,0}(Q)$, which is used in the following lemma.

\begin{lemma} \label{Welle:Lem:atilde_3Bed}
  For all $0 \neq v \in H^{1,1}_{0;,0}(Q)$, there exists a function
  $u_v \in \mathcal H_{0,}(Q)$ such that
  \[
    \widetilde{a}(u_v,v) > 0 \, .
  \]
\end{lemma}

\proof{For $0 \neq v \in H^{1,1}_{0;,0}(Q)$, there exists a unique solution
  $u_v \in H^{1,1}_{0;0,}(Q) \subset \mathcal H_{0,}(Q)$, satisfying
  \[
    a(u_v,w)=
    \langle v , w \rangle_{L^2(Q)} \quad
    \text{for all } w \in H^{1,1}_{0;,0}(Q).
  \]
  Using the representation \eqref{Welle:atilde_ist_a}, this gives
  \[
    \widetilde{a}(u_v,w) =
    \langle v , w \rangle_{L^2(Q)} \quad
    \text{for all } w \in H^{1,1}_{0;,0}(Q),
  \]
  and in particular for $w=v$, we obtain
  \[
    \widetilde{a}(u_v,v) = \| v \|_{L^2(Q)}^2 > 0,
  \]
  i.e., the assertion.}

\noindent
Next, we state the new variational setting for the wave equation
\eqref{Welle}. For given $f \in [H^{1,1}_{0;,0}(Q)]'$, we consider
the variational formulation to find $u \in \mathcal H_{0,}(Q)$ such that
\begin{equation}\label{Welle:VF_verallgemeinert}
  \widetilde{a}(u,v) =
  \langle f, v \rangle_Q  \quad \text{for all } v \in H^{1,1}_{0;,0}(Q),
\end{equation}
i.e., the operator equation
\begin{equation*}
  \mathcal E' \Box \widetilde{u} = f \quad \text{in }
  [H^{1,1}_{0;,0}(Q)]'.
\end{equation*}
With the properties of the bilinear form $\widetilde{a}(\cdot,\cdot)$, the
unique solvability of the variational formulation
\eqref{Welle:VF_verallgemeinert}, i.e., the main theorem of this paper,
is proven.

\begin{theorem} \label{Welle:Thm:ExistenzEindeutigkeit}
  For each given $f \in [H^{1,1}_{0;,0}(Q)]'$, there exists a unique solution
  $u \in \mathcal H_{0,}(Q)$ of the variational formulation
  \eqref{Welle:VF_verallgemeinert}. Furthermore,
  \begin{equation*}
    \mathcal L \colon \, [H^{1,1}_{0;,0}(Q)]' \to \mathcal H_{0,}(Q),
    \qquad \mathcal L f = u,
  \end{equation*}
  is an isomorphism satisfying 
  \begin{equation*}
    \| \Box \widetilde{u} \|_{[H^1_0(Q_-)]'} =
    \| \Box \widetilde{\mathcal Lf} \|_{[H^1_0(Q_-)]'} =
    \| f \|_{[H^{1,1}_{0;,0}(Q)]'}.
  \end{equation*}
\end{theorem}

\proof{With the help of the Banach--Ne\v{c}as--Babu\v{s}ka theorem
  \cite[Theorem 2.6]{ErnGuermond2004}, the results in
  \eqref{Welle:atilde_stetig}, \eqref{Welle:atilde_infsup} and
  Lemma~\ref{Welle:Lem:atilde_3Bed} yield the existence and uniqueness of
  the solution $u \in \mathcal H_{0,}(Q)$. In addition, with the variational
  formulation \eqref{Welle:VF_verallgemeinert}, the equalities
  \begin{equation*}
    \| \Box \widetilde{u} \|_{[H^1_0(Q_-)]'} =
    \sup\limits_{0 \neq v \in H^{1,1}_{0;,0}(Q)}
    \frac{| \widetilde{a}(u,v) |}{\| v \|_{H^{1,1}_{0;,0}(Q)}} =
    \sup\limits_{0 \neq v \in H^{1,1}_{0;,0}(Q)}
    \frac{| \langle f, v \rangle_Q |}{\| v \|_{H^{1,1}_{0;,0}(Q)}} =
    \| f \|_{[H^{1,1}_{0;,0}(Q)]'} .
\end{equation*}
hold true and, therefore, the assertion.}

\noindent
While for the ordinary differential equation \eqref{ODE}, this new approach
leads to the same variational setting as already considered in
\cite[Section 4]{SteinbachZank2020} and \cite{SteinbachZankFEM2019},
see Lemma~\ref{ODE:Lem:HNeu_ist_H10} and Corollary~\ref{ODE:Kor:VF_gleich},
the situation is different for the wave equation \eqref{Welle}. In greater
detail, for the variational formulation \eqref{Welle:VF_verallgemeinert}, the
Banach--Ne\v{c}as--Babu\v{s}ka theorem \cite[Theorem 2.6]{ErnGuermond2004}
is applicable, whereas the variational formulation \eqref{Welle:VF} does
not fit in this framework, see Theorem~\ref{Welle:Thm:InfSupH1}.
Additionally, Lemma~\ref{Welle:Lem:Darstellung_Bilinearform} and
\eqref{Welle:atilde_ist_a} show that the new variational formulation
\eqref{Welle:VF_verallgemeinert} is a generalization of the variational
formulation \eqref{Welle:VF}. Next, the following functions are given to
get a first impression of the solution space $\mathcal H_{0,}(Q)$.

\begin{remark}
  For $u \in C^2(\overline{Q})$ with $u_{|\Omega \times \{0\}}=u_{|\Sigma}=0$
  there holds $u \in H^{1,1}_{0;0,}(Q) \subset \mathcal H_{0,}(Q)$.
  Note that the second initial condition
  \begin{equation*}
    \partial_t u(\cdot,t)_{|t=0}  = 0 \quad \text{in } \Omega
  \end{equation*}
  is not incorporated in the ansatz space $\mathcal H_{0,}(Q)$,
  see \eqref{Welle:fv0}.
\end{remark}

\begin{remark}
  Consider the smooth function
  \begin{equation*}
    u(x,t) = \sin(\pi x) \sin(\pi t) \quad
    \text{ for } (x,t) \in (0,1)\times (0,1) = Q,
  \end{equation*}
  satisfying $u_{|\Omega \times \{0\}}=u_{|\Sigma} = 0$ and
  $\Box u = 0$ in $Q$. But there is
  $\Box \widetilde{u} \neq 0$ with
  \begin{equation*}
    \widetilde{u}(x,t) =
    \begin{cases}
      u(x,t) & \text{for } (x,t) \in Q, \\
      0 & \text{for } (x,t) \in Q_- \setminus Q,
    \end{cases}
  \end{equation*}
  since the distributional derivative fulfills
  \begin{equation*}
    \langle \Box \widetilde{u}, \varphi \rangle_{Q_-} =
    \int_0^T \int_\Omega u(x,t) \, \Box \varphi(x,t) \,
    \mathrm dx \, \mathrm dt =
    \pi \int_\Omega \sin(\pi x) \, \varphi(x,0) \, \mathrm dx
  \end{equation*}
  for all $\varphi \in C^\infty_0(Q_{-})$. Thus, the function
  $u \in H^{1,1}_{0;0,}(Q) \subset \mathcal H_{0,}(Q)$ solves the
  variational formulation \eqref{Welle:VF_verallgemeinert} with the
  right--hand side $f_{v_0} \in [H^{1,1}_{0;,0}(Q)]'$,
  \begin{equation*}
    \langle f_{v_0}, v \rangle_Q =
    \pi \int_\Omega \sin(\pi x) \, v(x,0) \, \mathrm dx, \quad
    v \in H^{1,1}_{0;,0}(Q),
  \end{equation*}
  i.e., the function $u$ satisfies the inhomogeneous initial condition
  \begin{equation*}
    \partial_t u(x,t)_{|t=0}  =  v_0(x) := \pi \sin(\pi x), \quad x \in \Omega,
  \end{equation*}
  see \eqref{Welle:fv0}.
\end{remark}

\section{Conclusions and outlook} \label{Sec:Zum}
In this paper, we presented a new approach to set up a bijection 
for the solution of the wave equation, when the right--hand side is
considered in the dual space of the test space of the variational
formulation. For this, we had to enlarge the ansatz space to prove
a related inf--sup stability condition. Based on these results, we
aim to derive a space--time finite element method for the numerical
solution of the wave equation, and of related problems, which is
unconditionally stable, and which also allows for an adaptive
resolution of the solution simultaneously in space and time, and for an
efficient solution, which is also parallel in time. First numerical
results are very promising, see \cite{DD26}, and the related numerical
analysis is ongoing work, and will be published elsewhere.

The presented results on the existence and uniqueness of solutions for
the wave equation, in particular the bijectivity results for the solution
operator in related function spaces, are of utmost importance for the
analysis of related boundary integral equations for the approximate
solution of the wave equation by boundary element methods. Using the
appropriate Dirichlet and Neumann trace operators, we are able to analyze
the mapping properties of related boundary integral operators
\cite{SteinbachUrzua2020}, i.e.,
boundedness and coercivity, to close the existing gap in using
different norms, see, e.g., \cite{Sayas2016}. Note that this norm gap
also results in error estimates, which are not optimal, see
also \cite{SteinbachUrzuaZank2020} for first numerical results.

We end this paper with an outlook for possible extensions of the approach
in Section~\ref{Sec:Welle}. Since the constructions of the spaces
$\mathcal H(Q)$, $\mathcal H_{0,}(Q)$ and the proofs in this section
mainly rely on the treatment of the second--order temporal differential
operator $\partial_{tt} + \mu$ with a parameter $\mu$, a generalization of
the results of this section to differential operators
$\partial_{tt} + \mathcal A_x,$ acting on vector fields or scalar
fields is possible, where the second--order spatial differential
operator $\mathcal A_x$ has to fulfill certain properties, e.g.,
boundedness and ellipticity. A more detailed discussion is left for
future work.


\end{document}